\documentclass[12pt]{amsart}

\usepackage{hyperref}
\hypersetup{nesting=true,debug=true,naturalnames=true}
\usepackage{graphicx,amssymb,upref,units}

\usepackage{tikz}

\let\<\langle
\let\>\rangle

\let\uml\"

\usepackage{enumitem,kantlipsum}

\usepackage{icomma}

\usepackage{amsmath,amssymb,fullpage,xy,mathrsfs,amsthm,amsxtra,graphicx,verbatim}
\usepackage{hyperref}
\hypersetup{colorlinks=true,urlcolor=magenta,citecolor=magenta,linkcolor=blue}

\usepackage{colonequals}
\usepackage{multirow}

\usepackage[mathcal]{euscript}
\newcommand{\Gal}{\operatorname{Gal}}
\newcommand{\End}{\operatorname{End}}

\newcommand{\OK}{\mathcal{O}_{K}}

\newcommand{\Fq}{\mathbf{F}_{q}}
\newcommand{\Aut}{\operatorname{Aut}}

\newcommand{\Z}{\mathbf{Z}}
\newcommand{\F}{\mathbf{F}}

\newcommand{\Q}{\mathbf{Q}}

\newcommand{\GL}{\operatorname{GL}}

\newcommand{\ddef}{\colonequals}
\newcommand{\disc}{\operatorname{disc}}

\newcommand{\tors}{\mathsf{tors}}

\newcommand{\sqf}{\mathsf{sqf}}
\newcommand{\head}{\mathrm{head}}
\newcommand{\tail}{\mathrm{tail}}
\newcommand{\PPP}{\mathcal{P}}
\newcommand{\Prop}{\operatorname{Prop}}

\numberwithin{equation}{section}

\theoremstyle{plain}
\newtheorem{thm}[equation]{Theorem}

\newtheorem{lem}[equation]{Lemma}
\newtheorem{defn}[equation]{Definition}
\newtheorem{cor}[equation]{Corollary}
\newtheorem{prop}[equation]{Proposition}

\theoremstyle{remark}
\newtheorem{rmk}[equation]{Remark}

\newtheorem{exm}[equation]{Example}

 \input xy
 \xyoption{all}

\begin{document}

\title{The probability of non-isomorphic group structures of isogenous elliptic curves in finite field extensions, II}

\author{John Cullinan}
\address{Department of Mathematics, Bard College, Annandale-On-Hudson, NY 12504, USA}
\email{cullinan@bard.edu}
\urladdr{\url{http://faculty.bard.edu/cullinan/}}

\author{Shanna Dobson}
\address{Department of Mathematics, University of California, Riverside, CA 92521, USA}
\email{Shanna.Dobson@email.ucr.edu}

\author{Linda Frey}
\address{Mathematisches Institut, Georg-August-Universit\"at G\"ottingen, G\"ottingen, Germany}
\email{lindafrey89@gmail.com}

\author{Asimina Hamakiotes}
\address{Department of Mathematics, University of Connecticut, Storrs, CT 06269, USA}
\email{asimina.hamakiotes@uconn.edu}

\author{Roberto Hernandez}
\address{Department of Mathematics, Emory University, Atlanta, GA 30322, USA}
\email{roberto.hernandez@emory.edu}

\author{Nathan Kaplan}
\address{Department of Mathematics, University of California, Irvine, CA 92697, USA}
\email{nckaplan@math.uci.edu}
\urladdr{\url{https://www.math.uci.edu/~nckaplan/}}

\author{Jorge Mello}	
\address{Department of Mathematics, Oakland University, Rochester , MI 48309, USA}
\email{jorgedemellojr@oakland.edu}

\author{Gabrielle Scullard}
\address{Department of Mathematics, The Pennsylvania State University, University Park, Pennsylvania 16802, USA}
\email{gns49@psu.edu}
\urladdr{\url{https://sites.psu.edu/gns49/}}

\begin{abstract}
Let $E$ and $E'$ be 2-isogenous elliptic curves over $\Q$.  Following \cite{ck}, we call a good prime $p$ \emph{anomalous} if $E(\F_p) \simeq E'(\F_p)$ but  $E(\F_{p^2}) \not \simeq E'(\F_{p^2})$.  Our main result is an explicit formula for the proportion of anomalous primes for any such pair of elliptic curves.  We consider both the CM case and the non-CM case.  
\end{abstract}

\maketitle

\section{Introduction}

\subsection{Motivation} In this paper we follow up on our work in \cite{ck} by analyzing the possible values for the proportions of so-called anomalous primes (as introduced in \cite{ck}). Let $E$ and $E'$ be 2-isogenous elliptic curves defined over $\Q$.  We call a good prime $p>2$ for $E$ and $E'$ \emph{anomalous} if $E(\F_p) \simeq E'(\F_p)$ but $E(\F_{p^2}) \not \simeq E'(\F_{p^2})$.  

The proportion that we study in this paper is defined concretely as
\begin{align} \label{prop_def}
\PPP(E,E') =  \lim_{X\rightarrow \infty} \frac{\#\lbrace \text{anomalous}~p~\text{for $E$ and $E'$} \leq X \rbrace}{\pi(X)},
\end{align}
where $\pi$ is the prime-counting function.  That this limit exists is ultimately due to the Chebatorev Density Theorem, applied to the images $G$ and $G'$ of the 2-adic representations of $E$ and $E'$, respectively.  

Because $E$ and $E'$ are 2-isogenous over $\Q$, each has a rational 2-torsion point and so $G$ and $G'$ each have index at least $3$ in $\GL_2(\Z_2)$.  When $[\GL_2(\Z_2):G] = [\GL_2(\Z_2):G']=3$, we showed that the proportion $\mathcal{P}(E,E') = 1/30$ \cite[Thm.~1.3.4]{ck}.  The purpose of this paper is to extend that result to all pairs of rationally 2-isogenous elliptic curves over $\Q$ and to explore some of the arithmetic consequences.

We break the large task of working out $\PPP(E,E')$ for every pair $(E,E')$ of 2-isogenous elliptic curves over $\Q$ into two main steps; we set 
\begin{align} \label{head_tail}
\mathcal{P}(E,E') = \mathcal{P}(E,E')^\head + \mathcal{P}(E,E')^\tail
\end{align}
and evaluate each term separately. We delay a formal definition of $\mathcal{P}(E,E')^\head$ and  $\mathcal{P}(E,E')^\tail$ until Section \ref{main_results}.

\subsection{Background and First Definitions} In order to efficiently state our main results, we start by recalling some of the terminology of \cite{ck}.  Since $E$ and $E'$ are 2-isogenous, the only way we can have $E(\F_{p^2}) \not \simeq E'(\F_{p^2})$ is if their 2-Sylow subgroups are non-isomorphic.  A convenient way to organize the anomalous primes is by their so-called \emph{defect}, introduced in \cite{ck}, which we now recall.

\begin{defn} \label{defect_defn}
We say that an anomalous prime $p$ has  \emph{defect} $(a,b)$ if $a \ne b$, $E(\F_{p^2})$ has full $2^a$-torsion but not full $2^{a+1}$-torsion, and $E(\F_{p^2})$ has full $2^b$-torsion, but not full $2^{b+1}$-torsion.
\end{defn}

We proved in \cite[Thm.~3.0.10]{ck} that the defect of a prime can only be of the form $(m+1,m)$ or $(m,m+1)$, where $m \geq 2$.  We organize primes by defect via the 2-adic representation.  Let $F$ and $F'$ be matrices that represent the Frobenius endomorphism on $E$ and $E'$, respectively, as elements of $\GL_2(\Z_2)$. 

\begin{prop} \label{defect_characterization}
A good prime $p$ for $E$ and $E'$ has defect $(m+1,m)$ if and only if 
\begin{align*}
&F \equiv -I \pmod{2^m} \text{ and } F \not \equiv -I \pmod{2^{m+1}}, \text{ and}\\
&F'\equiv -I \pmod{2^{m-1}} \text{ and } F' \not \equiv -I \pmod{2^{m}},
\end{align*}
with similar congruences for primes of defect $(m,m+1)$.
\end{prop}

\begin{rmk}
Proposition \ref{defect_characterization} was proved in \cite[\S5.1]{ck} with the stated hypothesis that $G$ and $G'$ have index 3 in $\GL_{2}(\Z_2)$.  In fact, this hypothesis is not used in the proof. 
\end{rmk}

\begin{defn}
Let $k \geq 2$.  We define $c_{2^k} \in [0,1]$ to be the proportion of primes such that $F' \not \equiv -I \pmod{2^k}$ given that $F \equiv -I \pmod{2^k}$.  Similarly, we define $c_{2^k}'$ to be the proportion of primes such that  $F \not \equiv -I \pmod{2^k}$, given that  $F' \equiv  -I \pmod{2^k}$. 
\end{defn}

Let $G$ and $G'$ denote the respective images of the 2-adic representations of $E$ and $E'$, and $G(2^k)$ and $G'(2^k)$ the respective images of the mod $2^k$ representations.  When $G$ and $G'$ have index 3 in $\GL_2(\Z_2)$, we determined that $\mathcal{P}(E,E') = 1/30$ by proving that $c_{2^k} = c_{2^k}' = 1/2$ for all $k \geq 2$ and then expressing $\mathcal{P}(E,E')$ as the sum of Chebotarev densities
\begin{align} \label{gen_form}
\PPP(E,E') = \sum_{k \geq 2} \left(\frac{c_{2^k}}{|G(2^k)|} + \frac{c'_{2^k}}{|G'(2^k)|}\right).
\end{align}
If $G$ and $G'$ have index 3 in $\GL_2(\Z_2)$, then $|G(2^k)| = |G'(2^k)| = 2^{4k-3}$, from which the value of 1/30 follows immediately by summing a geometric series.   By Proposition \ref{defect_characterization}, the formula (\ref{gen_form}) holds for general $G$ and $G'$, and it is this series that we study throughout the paper. 

\subsection{Main Results} \label{main_results} With our terminology in place, we can now turn to the main business of the paper.  Our first result is on the constants $c_{2^k}$ and $c_{2^k}'$.  

\begin{prop} \label{const_prop}
For all $k\geq 2$, the constants $c_{2^k}$ and $c_{2^k}'$ of \eqref{gen_form} satisfy $c_{2^k},c_{2^k}' \in \lbrace 0,1/2,1 \rbrace$ and $0 \leq c_{2^k} + c_{2^k}' \leq 1$.
\end{prop}

Roughly, the proof amounts to re-interpreting the conditional probabilities $c_{2^k}$ and $c_{2^k}'$ in terms of certain field extensions.  The fact that their values are restricted to $0$, $1/2$, and $1$ and sum is bounded by 1 relies on the fact that $|G(2^k)|/|G'(2^k)| \in \lbrace 2, 2^{-1},1 \rbrace$, which ultimately is due to the curves being 2-isogenous. 

Turning to the calculation of (\ref{gen_form}) we recall the fact that if $E/\Q$ is an elliptic curve and $-I \in G(32)$, then $-I \in G$ \cite[Cor.~1.3]{rzb}. This allows us to split $\PPP(E,E')$ into two terms that we label the `head' and `tail' as in (\ref{head_tail}).  Specifically, we define
\begin{align*}
\PPP(E,E')^\head  &= \sum_{k=2}^4 \left(\frac{c_{2^k}}{|G(2^k)|} + \frac{c_{2^k}'}{|G'(2^k)|}\right) \\
\PPP(E,E')^\tail  &= \sum_{k\geq 5} \left(\frac{c_{2^k}}{|G(2^k)|} + \frac{c_{2^k}'}{|G'(2^k)|}\right),
\end{align*}
and prove the following.

\begin{thm} \label{rough_thm}
Let $E$ and $E'$ be 2-isogenous elliptic curves defined over $\Q$ with 2-adic Galois images $G$ and $G'$, respectively.  Suppose $E$ and $E'$ do not have CM.  Then  $\PPP(E,E')^\tail =0$ if and only if $-I \not \in G$.  If $-I \in G$, then 
\[
\mathcal{P}(E,E')^\tail = \frac{8}{15}  \left(\frac{1}{|G(32)|} + \frac{1}{|G'(32)|} \right).
\]
\end{thm}

The proof of Theorem \ref{rough_thm} comes down to showing that either
\begin{itemize}
\item for every $k \ge 5$ we have $c_{2^k} = c'_{2^k} = 0$ (which occurs if and only if $-I \not \in G$), or
\item for every $k \ge 5$ we have $c_{2^k} = c'_{2^k} = 1/2$ (which occurs if and only if  $-I \in G$).
\end{itemize}
This observation relies on the fact that for elliptic curves over $\Q$, the largest 2-power level of a Galois representation is 32.  

\begin{rmk}
Theorem \ref{rough_thm} gives an easy procedure for computing $\PPP(E,E')$ for any pair of 2-isogenous elliptic curves: first check if $-I \in G$, then $\PPP(E,E')^\head$ is determined by calculating $c_{2^k},c_{2^k}'$ for $k=2,3,4$, which can be done by machine computation.
\end{rmk}

If an elliptic curve over $\Q$ has CM, then so does any $\Q$-isogenous elliptic curve.  Thus, when considering the ``CM case'' we can safely assume that both $E$ and $E'$ have CM.  If $E$ and $E'$ have CM, then the formulas for $\PPP(E,E')^{\rm head}$ and $\PPP(E,E')^{\rm tail}$ are significantly simplified and there are only two possible values of $\mathcal{P}(E,E')$, as shown in the following theorem (by contrast, there are many possible values of $\PPP(E,E')$ in the non-CM case).

\begin{thm} \label{cm_thm}
Let $E$ and $E'$ be 2-isogenous elliptic curves over $\Q$ with CM.  Then $\PPP(E,E') = 0$ or  $1/12$. 
\end{thm}

\begin{rmk}
In \cite{alvaro}, the author proves that there are precisely 28 possible 2-adic images attached to elliptic curves with CM over $\Q$.  In Section \ref{1/12}, in the course of proving Theorem \ref{cm_thm}, we show exactly which pairs of 2-adic images $(G,G')$ yield the proportion 1/12.  All remaining pairs of images yield a proportion of 0.
\end{rmk}

\begin{rmk}
In all cases, if $\PPP(E,E') = 0$ then there are exactly zero anomalous primes, as opposed to a set of density 0.
\end{rmk}

\subsection{Notation} If $E$ is an elliptic curve over $\Q$ then we write 
\begin{align*}
&\rho_{E,\ell}: \Gal_\Q \to \Aut(T_\ell E),\text{ and} \\
&\overline{\rho}_{E,\ell^n}: \Gal_\Q \to \Aut(T_\ell E/\ell^nT_\ell E), 
\end{align*}
for the $\ell$-adic and mod $\ell^n$ representations of $E$, respectively.  We write $G$ for the image of $\rho_{E,\ell}$ and $G(\ell^n)$ for the image of $\overline{\rho}_{E,\ell^n}$.  In the context of isogenous curves $E$ and $E'$ we will write $G'$ and $G'(\ell^n)$ for the corresponding images attached to $E'$. 

Let $p$ be a good prime for $E$ and $E'$.  Then the Frobenius endomorphisms lift to a conjugacy class in characteristic 0.  We will often denote by $F$ and $F'$ the images of these endomorphisms in $G(2^m)$ and $G'(2^m)$, respectively; the value of $m$ will usually be evident from context. 

We rely heavily on the database of 2-adic images in \cite{rzb} and also the LMFDB \cite{lmfdb} for generating examples.  When referring to specific examples of elliptic curves or number fields, we use LMFDB notation.  We performed our computer experiments using \textsf{Magma} \cite{magma}.

\subsection{Acknowledgments} This project was started at the Rethinking Number Theory Workshop in Summer 2022; we would like to thank the organizers for a stimulating research environment and support.  We would also like to thank \'Alvaro Lozano-Robledo for helpful conversations and Jeremy Rouse for supplying some of the examples in Section \ref{extremal}.  Shanna Dobson is supported by an NSF Mathematical Sciences Graduate Research Fellowship.  Nathan Kaplan was supported by NSF Grants DMS 1802281 and DMS 2154223.

\section{Background}

\subsection{Prior Work}
We start with a brief contextualization of this project. Let $\ell$ be a prime number.  If $E$ and $E'$ are $\ell$-isogenous elliptic curves defined over a finite field $k$ of characteristic $p \ne \ell$, then one can ask whether or not $E(k) \simeq E'(k)$ implies $E(K) \simeq E'(K)$ as $K$ ranges over finite extensions of $k$.  In \cite{cull1,cull2} the author showed that if the kernel of the isogeny $E \to E'$ is generated by a $k$-rational $\ell$-torsion point, then
\begin{itemize}[wide, labelwidth=!, labelindent=0pt]  
\item if $\ell$ is odd then $E(k) \simeq E'(k)$ implies $E(K) \simeq E'(K)$ for all finite extensions $K/k$, and 
\item if $\ell =2$, then if $E(k) \simeq E'(k)$ and $E(\mathbf{k}) \simeq E'(\mathbf{k})$, where $\mathbf{k}$ is the unique quadratic extension of $k$, then  $E(K) \simeq E'(K)$ for all finite extensions $K/k$.
\end{itemize}
This explains our interest in elliptic curves where $\ell = 2$ and $E(k) \simeq E'(k)$ and $E(\mathbf{k}) \not \simeq E'(\mathbf{k})$.

Given the situation when $\ell =2$, one can ask \emph{how often} $E(k) \simeq E'(k)$ and $E(\mathbf{k}) \not \simeq E'(\mathbf{k})$.  In \cite{ck} we took the point of view of ``fix $E$, vary $p$'', by fixing a pair of 2-isogenous elliptic curves $E$ and $E'$ defined over $\Q$ for which the isogeny is also defined over $\Q$ (we call such a pair \textbf{rationally 2-isogenous}), and asked if the primes for which $E(\F_p) \simeq E'(\F_p)$ and $E(\F_{p^2}) \not \simeq E'(\F_{p^2})$ have a well-defined proportion.   Under the condition that $G$ and $G'$ are maximal in $\GL_2(\Z_2)$, we showed that this proportion exists and is computable.  It will follow from our work in this paper that $\mathcal{P}(E,E')$ always exists and is computable, due to the fact that the $c_{2^k}$ and $c_{2^k}'$ of Proposition \ref{const_prop} are given in terms of finite-degree field extensions.

As already noted, when $E$ and $E'$ do not have CM and $G$ and $G'$ are maximal, given the hypothesis that $E$ and $E'$ have a rational 2-torsion point, then this proportion is 1/30.  This value of 1/30 can be  suggestively written as
\begin{align*}
\frac{1}{30} &=  \left( \frac{1}{2}\cdot \frac{1}{32} + \frac{1}{2} \cdot \frac{1}{32\cdot 16} + \cdots +\frac{1}{2} \cdot \frac{1}{32\cdot 16^m} + \cdots \right) \\ & +  \left( \frac{1}{2}\cdot \frac{1}{32} + \frac{1}{2} \cdot \frac{1}{32\cdot 16} + \cdots + \frac{1}{2} \cdot \frac{1}{32\cdot 16^m} + \cdots \right).
\end{align*}
The parenthetical expressions are identical, owing to the fact that $G$ and $G'$ are isomorphic when they have index 3 in $\GL_2(\Z_2)$.  Within each expression the factors of 1/2 are the $c_{2^k}$ for $k \geq 2$, and the fact that $|G(2^k)| = |G'(2^k)| = 2^{4k-3}$ is due to the maximality of $G$ and $G'$.  The purpose now is to generalize this to all  pairs of rationally 2-isogenous elliptic curves over $\Q$ to understand the possible values of $\PPP(E,E')$ that can occur.

\subsection{Definitions, Terminology, and Setup} \label{background}

To start, we would like to understand how much of the above analysis carries over to arbitrary 2-isogenous $E$ and $E'$; fortunately, nearly all of it does.  We start by collecting some well-known properties of 2-adic Galois representations of elliptic curves over $\Q$.  Our treatment is brief and, because much of this is well-known, we do not aim to be encyclopedic.

If $E$ is an elliptic curve over $\Q$ without CM, then Serre's open image theorem \cite{serre_open} implies that $G$ is a finite index subgroup of $\GL_2(\Z_2)$.   If $E$ and $E'$ are isogenous elliptic curves, then the indexes $[\GL_2(\Z_2):G]$ and $[\GL_2(\Z_2):G']$ are equal \cite[Prop.~2.2.1]{greenberg}.  Recall if $\pi_n: \GL_2(\Z_2) \to \GL_2(\Z/2^n\Z)$ denotes the reduction modulo $2^n$ homomorphism, then the \textbf{level} of $G$ is the smallest value $2^l$ such that $G = \pi_{l}^{-1} G(2^l)$.  For elliptic curves over $\Q$ the index divides 192 and the level divides 32 \cite{rzb}.  Up to isomorphism, there are 1208 possibilities for $G$ attached to elliptic curves over $\Q$ \cite{rzb}.  If $E$ has CM over $\Q$, then $G$ does not have finite index in $\GL_2(\Z_2)$, and the possibilities for $G$ have been classified in \cite{alvaro}; up to isomorphism there are 28 possibilities for $G$.  

Turning now to anomalous primes, we recall some terminology and results from \cite{ck}.  If $p$ is anomalous then by \cite[Lem.~3.0.8]{ck} we have 
\begin{align} \label{22}
E(\F_p)[2^\infty] \simeq E'(\F_p)[2^\infty] \simeq \Z/2\Z \times \Z/2\Z.
\end{align}
This immediately gives us the following.

\begin{lem} \label{4lem}
With all notation as above, if either $E(\Q)$ or $E'(\Q)$ has a rational 4-torsion point, then there are no anomalous primes.  
\end{lem}

\begin{proof} 
If either $E$ or $E'$ has a rational 4-torsion point, then it has an $\F_p$-rational 4-torsion point for all good primes $p$. By \eqref{22}, 
no good prime can be anomalous.
\end{proof}

\begin{cor} \label{tor_cor}
If $E$ and $E'$ are rationally 2-isogenous and there exists an anomalous prime, then 
$\Z/2\Z \subseteq E(\Q)[2^\infty], E'(\Q)[2^\infty] \subseteq \Z/2\Z \times \Z/2\Z$.
\end{cor}

\begin{rmk}
Lemma  \ref{4lem} and Corollary \ref{tor_cor} hold for any rationally 2-isogenous elliptic curves over $\Q$, whether or not they have CM.
\end{rmk}

\subsection{Frobenius Properties of Anomalous Primes} If $E$ and $E'$ are $\ell$-isogenous elliptic curves defined over a finite field, then the index $[\End(E):\End(E')]$ can take one of three values: $\ell$, $\ell^{-1}$, or 1. We call the isogeny $\varphi: E \to E'$ ascending/descending/horizontal in these cases, respectively.  It is also the case that over $\Q$,  $\ell$-isogenous curves have $|G(\ell^m)|/|G'(\ell^m)| \in \lbrace 1,\ell,\ell^{-1} \rbrace$ \cite[Prop.~2.1.1]{greenberg}.

Now let $E$ and $E'$ be rationally 2-isogenous elliptic curves.  If $p$ is anomalous, then $\varphi$, when reduced modulo $p$, is either descending or ascending.   If $\varphi$ were horizontal then $\End(E) \simeq \End(E')$, which would imply $E(\F_{p^k}) \simeq E(\F_{p^k})$ for all $k \geq 1$, since for ordinary primes we have 
\begin{align}
E(\F_{p^k}) \simeq \End(E)/(\pi^k - 1),
\end{align}
where $\pi$ is the Frobenius endomorphism \cite[Thm.~1]{lenstra}.  Suppose 
\[
|E(\F_{p^2})[2^\infty]| = |E'(\F_{p^2})[2^\infty]| = 2^v.
\]
Consider the Frobenius endomorphisms $F$ and $F'$ of $E$ and $E'$ at $p$ as elements in $G$ and $G'$, respectively.  What follows for the remainder of this section is a recap from \cite{ck} of the characterization of Frobenius at anomalous primes.  

If $\varphi$ is descending and $p$ is anomalous, then we must have 
\begin{align*}
E(\F_{p^2})[2^\infty] &\simeq \Z/2^{m+1}\Z \times \Z/2^{v-m-1}\Z \\
E'(\F_{p^2})[2^\infty] &\simeq \Z/2^{m}\Z \times \Z/2^{v-m}\Z, 
\end{align*}
for some $m \geq 2$ and where $v \geq 2m+2$.  This setup occurs precisely when $F \equiv -I \pmod{2^m}$ (but $F \not \equiv -I \pmod{2^{m+1}}$) and $F' \equiv -I \pmod{2^{m-1}}$ (but $F' \not \equiv -I \pmod{2^m}$), by \cite[\S 5.1]{ck}.   Recall from Definition \ref{defect_defn} that we say $p$ has \emph{defect} $(m+1,m)$.  Similar statements hold for ascending isogenies by interchanging the roles of $E$ and $E'$.  Therefore, we determine the defect of primes by characterizing the conditions under which Proposition \ref{defect_characterization} holds.

\begin{rmk} \label{1mod2^m}
Observe that if $p$ is an anomalous prime of defect $(m+1,m)$ or $(m,m+1)$, then $p \equiv 1\pmod{2^m}$, since then $F \equiv -I \pmod{2^m}$, which has determinant congruent to $1 \pmod{2^m}$.
\end{rmk}

\subsection{Galois Theory of Torsion Fields} \label{tors_fields} Let $A$ be an elliptic curve defined over a field $k$.  Then $A$ can be defined explicitly by a Weierstrass equation.  For any positive integer $T$, let $x(A[T]) \subseteq \overline{k}$ denote the set of $x$-coordinates of the $T$-torsion subgroup of $A$.  Then the field $k(x(A[T]))$ is a Galois extension of $k$ and is  independent  of the choice of Weierstrass model for $A$.  Similarly, let $k(A[T])$ be the $T$-torsion field of $A$, obtained by adjoining to $k(x(A[T]))$ the $y$-coordinates of the $T$-torsion points of $A$; this field is also independent of the choice of Weierstrass model for $A$.

Turning to the situation at hand, fix $m \geq 2$ and let $E$ and $E'$ be 2-isogenous elliptic curves defined over $\Q$.  Let
\begin{align*}
K &= \Q(x(E[2^m])) \text{ and } L = \Q(E[2^m]) \\
K' &= \Q(x(E'[2^m])) \text{ and } L' = \Q(E'[2^m]).
\end{align*}
Then $K,K',L,L'$ are all Galois over $\Q$ and recall that composites of Galois extensions are Galois.  Recall further that the proportion of primes that split completely in a Galois extension is the inverse of the degree of the extension over $\Q$.

Consider the group $\langle \pm I \rangle \subseteq \GL_2(\Z_2)$ and its projection $\langle \pm I \rangle_{2^m}$ into $\GL_2(\Z/2^m\Z)$. Then $G(2^{m}) \cap \langle \pm I \rangle_{2^{m}} = \Gal(L/K)$, and we have 
\[
\xymatrix{
L \ar^{\Gal(L/\Q) \simeq G(2^m)}@/^3pc/@{-}[dd] \ar@{-}_{G(2^{m}) \cap \langle \pm I \rangle_{2^{m}}}[d] \\
K \ar@{-}[d]\\
\Q
}
\]
so that $[L:K] =2$ or 1, depending on whether $G(2^m)$ contains $-I \pmod{2^m}$ or not, respectively \cite[Fig.~5.5]{adelman}. 

We will use all of this as part of our setup to determine the constants $c_{2^m}$ and $c_{2^m}'$ of Theorem \ref{rough_thm}. In view of this, we  assume for now that $G(2^m)$ contains $-I \pmod{2^m}$.  (If not, then $c_{2^n} = 0$ for all $n \geq m$, by Proposition \ref{const_prop}).

Let $\lbrace P,Q\rbrace$ be a basis for $T_2E$ over $\overline{\Q}$, chosen such that $\lbrace P_{2^n}, Q_{2^n} \rbrace$ is a basis for $E[2^n]$ for all $n$, and such that $E' = E/\langle P_2 \rangle$, where $P_{2^n}$ and $Q_{2^n}$ are the reductions modulo $2^n$ of $P$ and $Q$, respectively.  We write $P',Q'$ for a compatible basis of $T_2E'$ such that the isogeny $\varphi:E \to E'$ maps $P$ to $P'$ and $Q$ to $Q'$.

If $p$ is a prime such that $F \equiv -I \pmod{2^n}$,  then $F(P_{2^n}) = -P_{2^n}$ and $F(Q_{2^n}) = -Q_{2^n}$. Suppose that $E \to E'$ is a descending isogeny at $p$ (if it is ascending, we can interchange $E$ and $E'$ via the dual isogeny).  Then we  have  $Q_{2^n}' = Q_{2^n} + \langle P_2 \rangle$ and so
\[
F'(Q_{2^n}') = F(Q_{2^n} + \langle P_2 \rangle) = -Q_{2^n} + \langle P_2 \rangle = -Q_{2^n}'.
\]
Therefore, in a suitable ordering of the bases, if $F \equiv -I \pmod{2^n}$, then $F' \equiv \left( \begin{smallmatrix} -1 & * \\ 0 & * \end{smallmatrix} \right) \pmod{2^n}$.   Computing the determinant of $F$ modulo $2^n$ shows that $p \equiv 1 \pmod{2^n}$.  This implies
\begin{align} \label{F'}
F' \equiv  \begin{pmatrix} -1 & * \\ 0 & -1 \end{pmatrix} \pmod{2^{n}}. 
\end{align}
Because $E$ has full $2^n$-torsion defined over $\F_{p^2}$ (by squaring $F$), it must be the case that $E'$ has full $2^{n-1}$-torsion defined over $\F_{p^2}$.  Hence $F' \equiv -I \pmod{2^{n-1}}$ and so $* \equiv 0 \pmod{2^{n-1}}$ by squaring matrices. 

If $F \equiv -I \pmod{2^n}$, then the polynomial generating the $x$-coordinates $x(E[2^n])$ of $E[2^n]$ splits modulo $p$ \cite[\S 5.3]{adelman}, hence $\F_p(x(E[2^n]))  = \F_p$.  On the other hand, since $F \not \equiv I \pmod{2^n}$ it is not the case that $\F_p(E[2^n]) = \F_p$.  Because the $2^n$-torsion field of an elliptic curve is at most a quadratic extension of a field containing $x(E[2^n])$, it follows that in our setup $[\F_p(E[2^n]):\F_p]=2$.

Since the $x$-coordinate of $Q_{2^n}$ is  defined over $\F_p$, so is the $x$-coordinate of $Q'_{2^n}$.  We have established that $F' \equiv -I \pmod{2^{n-1}}$, so it is also the case that $\F_p$ contains $x(E'[2^{n-1}])$.  Therefore, the obstruction to $x(E'[2^n])$ being defined over $\F_p$ is whether or not the preimages of $P'_{2^{n-1}}$ under duplication are all defined over $\F_p$.  

For convenience, we suppose $p \geq 5$ (omitting a finite set of primes will not affect our density computations).  Let $(\xi,\eta)$ be the coordinates of $P_{2^{n-1}}'$ and let $E'$ have Weierstrass equation 
\[
y^2 = x^3 + ax+b.
\]
Then the $x$-coordinates of the preimages of $P'_{2^{n-1}}$ under duplication are given by the solutions to 
\begin{align*} 
\frac{x^4-2ax^2-8bx+a^2}{4y^2} = \xi
\end{align*}
which are the same as the roots of the quartic
\begin{align} \label{miret}
x^4 - 4\xi x^3 - 2ax^2 + (-4\xi a - 8b)x + (a^2 - 4\xi b).
\end{align}
This polynomial is either irreducible, splits into two conjugate quadratics, or splits completely over $\F_p$, all of which have implications for anomalous primes.  See \cite[Lem.~5.1.8]{ck} for a proof of this statement and \cite[\S5]{ck} for a more general discussion of the implications.

\section{Prime Decompositions} \label{prime_decomp} 

The purpose of this section is to work out the possible values of the constants $c_{2^m}$ and $c_{2^m}'$ and determine the conditions under which they occur.  For the remainder of this section we fix an integer $m \geq 2$.  We begin with a lemma.

\begin{lem}
Suppose neither $G(2^m)$ nor $G'(2^m)$ contains $-I \pmod{2^m}$.  Then $c_{2^k} = c_{2^k}'=0$ for all $k \geq m$.
\end{lem}

\begin{proof}
Let $k \geq m$.  By Proposition \ref{const_prop}, a prime is anomalous of defect $(k+1,k)$ or $(k,k+1)$ only if $F \equiv -I \pmod{2^k}$ or $F' \equiv -I \pmod{2^k}$, respectively.  If neither $G(2^m)$ nor $G'(2^m)$ contains $-I \pmod{2^m}$, then neither $G(2^k)$ nor $G'(2^k)$ contains $-I \pmod{2^k}$.  Thus there can be no primes of defect $(k+1,k)$ or $(k,k+1)$ and so $c_{2^k} = 0$ and $c_{2^k}'=0$, accordingly.
\end{proof}

Next, we determine the values of $c_{2^m}$ and $c_{2^m}'$ when exactly one of $G(2^m)$ or $G'(2^m)$ contain $-I \pmod{2^m}$.

\begin{lem}
If $G(2^m)$ contains $-I \pmod{2^m}$ and $G'(2^m)$ does not, then $c_{2^m} = 1$ and $c_{2^k}' = 0$ for all $k \geq m$.  Similarly, if $G'(2^m)$ contains $-I \pmod{2^m}$ and $G(2^m)$ does not, then $c_{2^m}' = 1$ and $c_{2^k} = 0$ for all $k \geq m$.
\end{lem}

\begin{proof}
We prove the first claim and omit the identical argument for the second.  A prime is anomalous of defect $(k+1,k)$ only if $F' \equiv -I \pmod{2^k}$.  Since $G'(2^m)$ does not contain $-I \pmod{2^m}$, it is also true that $G'(2^k)$ does not contain $-I \pmod{2^k}$ for all $k \geq m$.   Thus there can be no primes of defect $(k,k+1)$, and so $c_{2^k}' =0$ for all $k \geq m$.  

On the other hand, for every prime $p$ such that $F \equiv -I \pmod{2^m}$, it can never be the case that $F' \equiv -I \pmod{2^m}$, hence every such prime is anomalous of defect $(m+1,m)$.  Therefore, $c_{2^m}=1$.
\end{proof}

For the remainder of this section we will assume that both $G(2^m)$ and $G'(2^m)$ contain $-I \pmod{2^m}$.  We observe from the start that this implies $[L:K] = [L':K']=2$ from our work in Section \ref{tors_fields}. The values of $c_{2^m}$ and $c_{2^m}'$ will be determined by the splittings of primes in certain composite field extensions, since each of those values depends on both $G(2^m)$ and $G'(2^m)$.  We start by setting up the structure of these composite fields.

A rational prime $p$ splits completely in $K$ if and only if $F \equiv \pm I \pmod{2^m}$.  Similarly, $p$ splits completely in $K'$ if and only if $F' \equiv \pm I \pmod{2^m}$.  By the Chebatorev Density Theorem, the proportions of these primes are $2/|G(2^m)|$ and $2/|G'(2^m)|$, respectively.  Thus $2/|G(2^m)|$ of the primes split completely in $K$, $2/|G'(2^m)|$ of the primes split completely in $K'$, and we verify that
\[
[K:\Q] = \frac{|G(2^m)|}{|\lbrace \pm I \rbrace|} = \frac{|G(2^m)|}{2} \text{\qquad and \qquad } [K':\Q] = \frac{|G'(2^m)|}{|\lbrace \pm I \rbrace|} = \frac{|G'(2^m)|}{2}.
\]
Similarly, a prime $p$ splits completely in $L$ if and only if $F \equiv I \pmod{2^m}$, and splits completely in $L'$ if and only if $F' \equiv I \pmod{2^m}$, and this occurs for $1/|G(2^m)|$ and $1/|G'(2^m)|$ of the primes, respectively.  

The proportion of primes for which $F \equiv I \pmod{2^m}$ and $F' \equiv I \pmod{2^m}$ is 
\begin{align*}
\frac{1}{[LL':\Q]} &= \frac{1}{[LL':L][L:\Q]} = \frac{1}{[LL':L]|G(2^m)|}\\ 
&= \frac{1}{[LL':L'][L':\Q]} = \frac{1}{[LL':L']|G'(2^m)|}.
\end{align*}
Similarly, the proportion of primes for which $F \equiv \pm I \pmod{2^m}$ and $F' \equiv \pm I \pmod{2^m}$ is 
\begin{align*}
\frac{1}{[KK':\Q]} &= \frac{1}{[KK':K][K:\Q]} = \frac{2}{[KK':K]|G(2^m)|}\\
&= \frac{1}{[KK':K'][K':\Q]} = \frac{2}{[KK':K']|G'(2^m)|}.
\end{align*}
Hence, 
\begin{align} \label{ratios}
\frac{[LL':L']}{[LL':L]} = \frac{[KK':K']}{[KK':K]}  = \frac{|G(2^m)|}{|G'(2^m)|} \in \lbrace 1,2,1/2 \rbrace,
\end{align}
where the last claim follows from the fact (referenced in Section \ref{background}) that $\ell$-isogenous curves have $|G(\ell^m)|/|G'(\ell^m)| \in \lbrace 1,\ell,\ell^{-1} \rbrace$.

\begin{lem} \label{sameF}
It cannot be the case that $F \equiv I \pmod{2^m}$ and $F' \equiv -I \pmod{2^m}$ or that $F \equiv -I \pmod{2^m}$ and $F' \equiv I \pmod{2^m}$.
\end{lem}

\begin{proof}
Since $E$ and $E'$ are $\Q$-isogenous, they are $\F_p$-isogenous when reduced modulo $p$ and so $|E(\F_p)| = |E'(\F_p)|$.   If $F \equiv I \pmod{2^m}$, then $v_2(|E(\F_p)|) \geq {2m}$, while if  $F' \equiv -I \pmod{2^{m}}$, then $v_{2}(|E'(\F_p)|) \leq 2$. Since $m \geq 2$, this is not possible (and neither is the case where $F \equiv -I \pmod{2^m}$ and $F' \equiv I \pmod{2^m}$).
\end{proof}

\begin{defn}
Let $X \in [0,1]$ be the proportion of primes $p$ such that $F' \equiv -I \pmod{2^m}$ given that $F \equiv -I \pmod{2^m}$.  Similarly, let $X' \in [0,1]$ be the proportion of primes $p$ such that $F \equiv -I \pmod{2^m}$ given that $F' \equiv -I \pmod{2^m}$.
\end{defn}

\begin{prop} \label{Xprop}
With all notation as above, we have
\begin{align} \label{Xeq}
X &= \frac{2}{[KK':K]} - \frac{1}{[LL':L]}, \text{ and}\\
X' &= \frac{2}{[KK':K']} - \frac{1}{[LL':L']}.
\end{align}
\end{prop}

\begin{proof}
Consider the proportion $X$. Writing $\Prop$ for the proportion of a set of primes, we have already seen 
\[
\Prop(F \equiv \pm I \pmod{2^m} \text{ and } F' \equiv \pm I \pmod{2^m}) = \frac{2}{[KK':K]|G(2^m)|}.
\]
In view of Lemma \ref{sameF}, we can rewrite this equation as the sum of conditional probabilities
\[
\Prop(F \equiv F' \equiv I \pmod{2^m}) + \Prop(F \equiv F' \equiv -I \pmod{2^m}) = \frac{2}{[KK':K]|G(2^m)|}.
\]
We can express this in terms of field extensions and the unknown quantity $X$:
\[
\frac{1}{[LL':L]|G(2^m)|} + \frac{X}{|G(2^m)|} = \frac{2}{[KK':K]|G(2^m)|},
\]
whence
\begin{align} \label{Xeq2}
X = \frac{2}{[KK':K]} - \frac{1}{[LL':L]}.
\end{align}
An identical argument gives the expression for $X'$.
\end{proof}

\begin{lem}
We have $c_{2^m} = 1-X$ and $c_{2^m}' = 1-X'$.
\end{lem}

\begin{proof}
By definition, $c_{2^m}$ is the proportion of primes $p$ with $F \equiv -I \pmod{2^m}$ such that $p$ is anomalous, and $X$ is the proportion of primes $p$ with $F \equiv -I \pmod{2^m}$ and $F' \equiv -I \pmod{2^m}$.  But $p$ is anomalous of defect $(m+1,m)$ if and only if $F \equiv -I \pmod{2^m}$ and $F' \not \equiv -I \pmod{2^m}$, hence $c_{2^m} = 1-X$.  An identical argument shows that $c_{2^m}' = 1-X'$.
\end{proof}

Each of the extensions $KK'/K$ and $LL'/L$ is generated by the same quartic \eqref{miret}, defined over $K$ (and hence over $L$), hence $[KK':K],[LL':L] \leq 4$. But because all Galois groups in the tower over $\Q$ are 2-groups, we can only have $[KK':K],[LL':L] \in \lbrace 1,2,4 \rbrace$.  Similarly, $[KK'\colon K'], [LL':L'] \in \lbrace 1,2,4 \rbrace$. This gives nine values each for $X$ and $X'$, six of which lie outside the range $[0,1]$, leaving us with the following twelve possibilities:
\begin{center}
\begin{tabular}{c|c|c||c|c|c}
$[KK':K]$ & $[LL':L]$ & $X$ & $[KK':K']$ & $[LL':L']$ & $X'$ \\
\hline
1 & 1 & 1 & 1 & 1 & 1\\
 2 & 1 & 0 & 2 & 1 & 0 \\
2 & 2 & 1/2 & 2 & 2 & 1/2 \\
4 & 2 & 0 & 4& 2& 0\\
2 & 4 & 3/4 &  2 & 4& 3/4\\
4 & 4 & 1/4 & 4&4&1/4
\end{tabular}
\end{center}
Over the next several lemmas we will rule out many of these possibilities for simple reasons and, from the remaining cases, deduce all possible values for $c_{2^m}$ and $c_{2^m}'$ and the conditions under which they occur. 

\begin{lem} \label{Xlem2}
It cannot be the case that $[LL':L] = 4$ or $[LL':L'] = 4$.
\end{lem}

\begin{proof}
Consider the quartic polynomial (\ref{miret}) defined over $K$, and hence over $L$.  The roots of this polynomial generate $KK'/K$ and $LL'/L$, so $[LL':L] = 1,2$, or 4.  By \cite[(5), p.~414]{miret}, the quartic splits over $L$ into two conjugate quadratics, hence $[LL':L] \ne 4$.  An identical argument shows that $[LL':L'] \ne 4$.
\end{proof}

\begin{prop}
It cannot be the case that $[KK':K]=4$ or that $[KK':K']=4$.
\end{prop}

\begin{proof}
We prove the first claim and the second follows by interchanging the roles of $K$ and $K'$.  If $F \equiv -I \pmod{2^m}$, then by \eqref{F'}, we have that $F' \equiv \left( \begin{smallmatrix} -1 & * \\ 0 & -1 \end{smallmatrix} \right) \pmod{2^m}$.  However, we can  deduce more about the unknown quantity `$*$' from the fact that $E$ and $E'$ are 2-isogenous.

Because $E$ and $E'$ are rationally 2-isogenous, they are 2-isogenous modulo $p$ for all good primes $p$.  Thus, $|E(\F_p)| = |E'(\F_p)|$ and $|E(\F_{p^2})| = |E'(\F_{p^2})|$, and the 2-Sylow subgroups are constrained by the 2-isogeny. In particular, since $F^2 \equiv I \pmod{2^{m+1}}$, $E$ has full $2^{m+1}$-torsion over $\F_{p^2}$.  Because $E$ and $E'$ are 2-isogenous, $E'$ must have full $2^m$-torsion over $\F_{p^2}$, which implies $F'^2 \equiv I \pmod{2^m}$.  It follows from squaring matrices modulo $2^{m+1}$ that $E'$ has full $2^{m+1}$-torsion over $\F_{p^2}$ if and only if $F'\equiv -I \pmod{2^m}$, and if $E'$ does not have full $2^{m+1}$-torsion over $\F_{p^2}$ then 
\[
F' \equiv \begin{pmatrix} -1 & 2^{m-1} \\ 0 & -1 \end{pmatrix} \pmod{2^m}.
\]

In turn, this leaves us with two possibilities when considering the quartic \eqref{miret} generating $K'$, viewed as a polynomial over $K$: 
\begin{itemize}
\item The quartic may split over $K$, in which case $F \equiv -I \pmod{2^m}$.  This occurs if and only if $F' \equiv -I \pmod{2^m}$, and thus $KK' = K$.
\item The quartic may factor into two conjugate, irreducible quadratics, in which case $[KK':K] = 2$ (this implies that for half of the primes $E'$ will have full $2^{m+1}$-torsion over $\F_{p^2}$ and for the complementary half $E'$ will only have full $2^m$-torsion over $\F_{p^2}$).
\end{itemize}
In neither case do we have $[KK':K] = 4$.
\end{proof}

\begin{lem}
It cannot be the case that $[KK':K] = 2$ and $[LL':L] =1$.  It also cannot be the case that $[KK':K'] = 2$ and $[LL':L']=1$. 
\end{lem}

\begin{proof}
Similar to the previous results, we prove the first claim and the second follows by an identical argument. By our standing assumption, both $G(2^m)$ and $G'(2^m)$ contain $-I$, which implies $[L:K] = [L':K']= 2$.  Now suppose $[KK':K] = 2$ and $[LL':L] = 1$.  If $[LL':L'] = 2$, then we would have $K' = L'$, contradicting $[L':K'] = 2$.  Thus $[LL':L']=1$, and so we have $LL' = L = L'$.  Then
\[
[LL':K] = [L:K] =2 =  [LL':KK'][KK':K].
\]
The hypothesis $[KK':K] = 2$ implies $LL' = KK'$.  Altogether these fields fit into the following diagram:
\[
\xymatrix{
L = L' = LL' = KK' \ar@{-}_2[d] \\
K}
\]

If $\mathcal{E}$ is an elliptic curve over $\Q$ with mod $2^m$ Galois image $\mathcal{G}(2^m) \subseteq \GL_2(\Z/2^m\Z)$, then by  \cite[Prop.~5.2.2, Fig.~5.5]{adelman} we have that  $\Q(x(\mathcal{E}[2^m]))$ is the unique index-2 subfield of $\Q(\mathcal{E}[2^m])$ corresponding to the fixed-field of the normal subgroup $\lbrace \pm I \rbrace$ of $\mathcal{G}(2^m)$.  We apply this observation to the setup of this lemma.

Since $L = L'$, we must have $G(2^m) \simeq G'(2^m)$.  This implies $K=K'$ as well, since each can be described as the unique index-2 subfield of $L$ corresponding to $\lbrace \pm I \rbrace$.  But then $KK' = K = K'$, contradicting the assumption that $[KK':K]=2$. Thus we cannot have $[KK':K] = 2$ and $[LL':L] = 1$ simultaneously.
\end{proof}

We are now poised to work out all possible values of $c_{2^m}$ and the conditions under which they occur.  The following lemma will apply to each of the cases in Proposition \ref{c2mprop} below, so we separate it out next.

\begin{lem} \label{containment_lem}
With all notation as above, if $K \subseteq K'$ then $c_{2^m}'=0$ and if $K' \subseteq K$ then $c_{2^m}=0$.
\end{lem}

\begin{proof}
If $K \subseteq K'$ then every rational prime that splits in $K'$ also splits in $K$.  
Applying Lemma \ref{sameF} shows that every time $F' \equiv -I \pmod{2^m}$ it must also be the case that $F \equiv -I \pmod{2^m}$ as well.  Therefore, there cannot exist primes of defect $(m,m+1)$ and so $c_{2^m}'=0$.  Similarly, if $K' \subseteq K$ then we must have $c_{2^m}=0$. 
\end{proof}

Four cases remain: $[KK':K] = [LL':L] =1$, $[KK':K'] = [LL':L']=1$, $[KK':K] = [LL':L] =2$, and $[KK':K']=[LL':L']=2$. Each of these cases can occur and yields coefficients $c_{2^m},c_{2^m}' \in \lbrace 0,1/2,1 \rbrace$.

\begin{prop} \label{c2mprop}
Suppose $G(2^m)$ and $G'(2^m)$ each contain $-I \pmod{2^m}$.  With all field notation as above, we have the following cases:
\begin{enumerate}
\item Suppose $[KK':K] = [LL':L] =1$. 
\begin{enumerate}
\item If $K' = K$  then $c_{2^m} = c_{2^m}'=0$.
\item If $K' \ne K$ then $c_{2^m} = 0$ and $c_{2^m}'=1/2$.
\end{enumerate}
\item Suppose $[KK':K'] = [LL':L']=1$.  
\begin{enumerate}
\item If $K' = K$  then $c_{2^m} = c_{2^m}'=0$.
\item If $K' \ne K$ then $c_{2^m}' = 0$ and $c_{2^m}=1/2$.
\end{enumerate}
\item Suppose $[KK':K] = [LL':L] =2$.  
\begin{enumerate}
\item If $K \subseteq K'$ then $c_{2^m} = 1/2$ and $c_{2^m}' = 0$.
\item If $K \not \subseteq K'$ then $c_{2^m} = c_{2^m}' = 1/2$.
\end{enumerate}
\item Suppose $[KK':K'] = [LL':L'] =2$.  
\begin{enumerate}
\item If $K' \subseteq K$ then $c_{2^m}' = 1/2$ and $c_{2^m} = 0$.
\item If $K' \not \subseteq K$ then $c_{2^m}' = c_{2^m} = 1/2$.
\end{enumerate}
\end{enumerate}
\end{prop}

\begin{proof}
We prove cases (1) and (3), since (2) and (4) follow by identical arguments.  Recall our standing assumption that $[L\colon K] = [L'\colon K'] = 2$; it will be used throughout this proof.

\medskip

Suppose $[KK':K] = [LL':L]=1$.  Then we have $K' \subseteq K$, hence $c_{2^m}=0$ by Lemma \ref{containment_lem}.  If, in addition, $K \subseteq K'$ so that $K=K'$, then $c_{2^m}'=0$.  This proves 1(a).  Now consider the case of 1(b) and suppose $K \ne K'$.  Then $[KK':K']=2$, hence 
\[
X' = \frac{2}{2} - \frac{1}{[LL':L']}.
\]
If we can show $[LL':L']=2$, we are done.  Since
\[
[LL'\colon K] = [LL' \colon KK'] [KK'\colon K] = [LL'\colon L] [L\colon K] = 2, 
\]
we have $2 = [LL' \colon KK']$.  But this gives
\[
[LL'\colon K'] = [LL'\colon L'] [L'\colon K'] = [LL' \colon KK'] [KK'\colon K'] = 4.
\]
 So $[L'\colon K'] = 2$ implies $[LL'\colon L'] = 2$, proving 1(b).

\medskip

Now we turn to case 3.  Since $[KK':K] = [LL':L] = 2$, we immediately have $c_{2^m} = 1/2$ by Proposition \ref{Xprop}. Thus 
\[
[LL':K] = [LL':L][L:K] = 4,
\]
from which it follows that $[LL':KK'] = 2$.  Now we consider possibilities 3(a) and 3(b) separately: either $K \subseteq K'$ or $K \not \subseteq K'$.  

If $K \subseteq K'$, then  $c_{2^m}'=0$ by Lemma \ref{containment_lem}, proving 3(a).  If $K \not \subseteq K'$ then $KK' \ne K'$, whence $[KK':K']=2$.  Together with $[LL':KK']=2$, we see that $[LL':L']=2$ and so $X' = 2/2 - 1/2 = 1/2 = c_{2^m}'$, proving 3(b).
\end{proof}

Our next proposition shows that if the proportion of primes of defect $(m+1,m)$ or $(m,m+1)$ is zero, then there are exactly zero primes of defect $(m+1,m)$ or $(m,m+1)$, respectively.  We use this result in the proof of Theorem \ref{cm_thm}, where we show that if there exists a single prime of defect $(m+1,m)$, then there must be a positive proportion.

\begin{prop} \label{literal_zero}
Suppose there exists a prime of defect $(m+1,m)$ or $(m,m+1)$.  Then the proportion of primes of defect $(m+1,m)$ or $(m,m+1)$, respectively, is positive.
\end{prop}

\begin{proof}
If there exists a prime of defect $(m+1,m)$, then in the group $G(2^m) \times G'(2^m)$ the Frobenius at such a  prime is represented by the pair $(-I,F')$, where $F'$ is conjugate in $G(2^m)$ to $\left(\begin{smallmatrix} -1 & 2^{m-1} \\ 0 & -1 \end{smallmatrix} \right)$.  Let $\mathcal{C}$ denote the conjugacy class of $G(2^m) \times G'(2^m)$ to which $(-I, F')$ belongs.  Observe that conjugation on each component preserves the fact that the the first component acts as $-I$ and the second does not.

The Galois group $\Gal(LL'/\Q)$ is a subgroup of $G(2^m) \times G'(2^m)$ and, because there exists a prime of defect $(m+1,m)$, intersects $\mathcal{C}$ nontrivially.   Thus the class in $\Gal(LL'/\Q)$ containing $(-I,F')$ is nonempty and by the Chebatorev Density Theorem the primes of defect $(m+1,m)$ will have a positive density.
\end{proof}

\begin{rmk}
Observe that $\Gal(LL'/\Q)$ is typically (much) smaller than $G(2^m) \times G'(2^m)$.  For example, let $E$ be the elliptic curve \href{https://www.lmfdb.org/EllipticCurve/Q/63/a/2}{63.a2} and $E'$ the curve \href{https://www.lmfdb.org/EllipticCurve/Q/63/a/5}{63.a5}.  The mod 4 images of these curves have size 4 and 8, respectively.  A computation in \textsf{Magma} shows that the composite $LL'$ has degree 8 over $\Q$ with $L \subseteq L'$.
\end{rmk}

To conclude this section, we give examples to show that  all these cases can actually occur over $\Q$. For cases 3(b) and 4(b), we refer the reader to the proof of \cite[Thm.~5.2.1]{ck} (which shows that $\PPP(E,E') = 1/30$ when $[\GL_2(\Z_2):G] = [\GL_2(\Z_2):G']=3$). There it is shown that whenever $[\GL_2(\Z_2):G] = [\GL_2(\Z_2):G']=3$ we have $[KK':K] = [LL':L] =2$, $K \not \subseteq K'$, and $K' \not \subseteq K$. 

\begin{exm}
Let $E$ be the elliptic curve \href{https://www.lmfdb.org/EllipticCurve/Q/277440/dv/5}{277440.dv5} and $E'$ the curve \href{https://www.lmfdb.org/EllipticCurve/Q/277440/dv/6}{277440.dv6}.   Set $m=2$.  Then $K \subsetneq K'$ and $[KK':K]=[LL':L]=2$, which is Case 3(a) (considering the dual isogeny yields 4(a)). It is also the case that 
$[KK':K']=[LL':L']=1$ which serves as an example for  Case 2(b) (taking the dual yields 4(b)).
\end{exm}

\begin{exm}
Let $E$ be the elliptic curve \href{https://www.lmfdb.org/EllipticCurve/Q/2304/i/1}{2304.i1} and $E'$ the curve \href{https://www.lmfdb.org/EllipticCurve/Q/2304/i/2}{2304.i2}. These curves have CM and both $G$ and $G'$ are maximal, according to the classification of 2-adic CM images in \cite{alvaro}.  

Now set $m=4$ and consider the fields $K$, $K'$, $L$, and $L'$.  We check that 
\[
\Q(E[2]) \simeq \Q(E'[2]) \simeq \Q(\sqrt{2}),
\]
while $K$ is the splitting field over $\Q$ of the polynomial $x^4 + 240x^2 - 3584x + 14400$ and $K'$ the splitting field of $x^4 + 60x^2 + 448x + 900$; one checks that $K \simeq K'$.  Finally, $L$ and $L'$ are each isomorphic to the degree-16 field 
\href{https://www.lmfdb.org/NumberField/16.0.118192468620711297024.15}{16.0.118192468620711297024.15}.    Thus, $[LL':L] = [LL':L'] = [KK':K] = [KK':K'] = 1$ and $c_4 = c_4' =0$.  This serves as an example for Cases 1(a) and 2(a) for Proposition \ref{c2mprop} above.
\end{exm}

\section{Elliptic Curves without CM}

We assume throughout this section that $E$ and $E'$ do not have CM so that $G$ and $G'$ are finite index subgroups of $\GL_2(\Z_2)$ of level at most 32.  Recall that we write 
\[
\PPP(E,E') = \PPP(E,E')^{\rm head} + \PPP(E,E')^{\rm tail}.
\]
One of the key distinctions in \cite{rzb} on the classification of 2-adic images of elliptic curves over $\Q$ is between those groups that contain $-I$ and those that do not.  The next lemma shows that this is an important distinction for us as well. 

\begin{lem}
Suppose $G$ and $G'$ do not contain $-I$.  Then $\PPP(E,E')^{\rm tail}=0$. 
\end{lem}

\begin{proof}
If $-I \not \in G$ (resp.~$G'$), then $-I \not \in G(32)$ (resp.~$G'(32)$) and therefore $-I \not \in G(2^m)$ (resp.~$G'(2^m)$) for $m \geq 5$ .  By our review in Section \ref{background}, we conclude there can be no anomalous primes of defect $(m+1,m)$ or $(m,m+1)$ for $m\geq 5$, and so $\PPP(E,E')^{\rm tail} =0$.
\end{proof}

On the other hand, if $-I$ is in both $G$ and $G'$ then we are guaranteed the existence of primes of defect $(m+1,m)$ and $(m,m+1)$ for $m \geq 5$, as the next proposition shows.

\begin{prop} \label{setup_for_tail}
Suppose $G$ and $G'$ contain $-I$.  Then there exist primes of defect $(m,m+1)$ and $(m+1,m)$ for all $m \geq 5$.  
\end{prop}

\begin{proof}
Since $G$ and $G'$ contain $-I$, we see that $-I \in G(2^m)$ for all $m$.  The image of the 2-adic representation of an elliptic curve over $\Q$ without CM has level at most $32$, so for each $m \ge 5$, $G(2^{m+1})$ contains an element $F$ of the form
\[
F = -I + 2^m \begin{pmatrix} 1 & 1 \\ 1 & 0 \end{pmatrix}.
\]
By \cite[Thm.~5.2.4]{ck}, the existence of such matrices implies the existence of primes of defect $(m+1,m)$.  An identical argument applied to $G'$ shows there exist primes of defect $(m,m+1)$ as well. 
\end{proof}

\begin{thm} \label{tailthm}
Suppose $G$ and $G'$ contain $-I$.  Then 
\[
\PPP(E,E')^{\rm tail}  = \frac{8}{15} \left( \frac{1}{|G(32)|} + \frac{1}{|G'(32)|} \right).
\]
\end{thm}

\begin{proof}
Fix $m \geq 5$.  Since both $G$ and $G'$ contain $-I$, Proposition \ref{setup_for_tail} implies that there exist primes of defect $(m+1,m)$ and $(m,m+1)$.  Therefore, both $c_{2^m}$ and $c_{2^m}'$ are nonzero.  By Proposition \ref{c2mprop}, this means $c_{2^m} = c_{2^m}' = 1/2$.  

Consider the primes of defect $(m+1,m)$.   The proportion of primes for which $F \equiv -I \pmod{2^m}$ is $1/|G(2^m)|$ by the Chebatorev Density Theorem.  Since $G$ has level at most $32$, we have $|G(2^m)| = |G(32)|\cdot 16^{m-5}$.  Thus, the proportion of primes of defect $(m+1,m)$ is 
\[
\frac{1}{2} \cdot \frac{1}{|G(2^m)|} = \frac{1}{2} \cdot \frac{1}{|G(32)|\cdot 16^{m-5}}.
\]
Similarly, the proportion of primes of defect $(m,m+1)$ is $\left( 2 \cdot |G'(32)| \cdot 16^{m-5} \right)^{-1}$. Summing over $m \geq 5$ gives
\[
\PPP(E,E')^{\rm tail} = \sum_{m \geq 5}\left( \frac{1}{2} \cdot \frac{1}{|G(32)|} \cdot \frac{1}{16^{m-5}} +  \frac{1}{2} \cdot \frac{1}{|G'(32)|} \cdot \frac{1}{16^{m-5}}  \right) = \frac{8}{15} \left( \frac{1}{|G(32)|} + \frac{1}{|G'(32)|} \right).
\]
\end{proof}

It remains to determine $\mathcal{P}(E,E')^{\rm head}$, which is just a matter of putting together results already established.

\begin{thm} \label{headthm}
Fix $m$ satisfying $2 \leq m \leq 4$.  The proportion of primes of defect $(m+1,m)$ is $c_{2^m}/|G(2^m)|$ with $c_{2^m} \in \lbrace 0,1/2,1 \rbrace$.  Similarly, the proportion of primes of defect $(m,m+1)$ is $c_{2^m}'/|G'(2^m)|$ with $c_{2^m}' \in \lbrace 0,1/2,1 \rbrace$.
\end{thm}

\begin{proof}
This follows immediately from Proposition \ref{c2mprop}.
\end{proof}

\noindent Together, Theorems \ref{headthm} and \ref{tailthm} prove Theorem \ref{rough_thm}.

\section{Extremal Values} \label{extremal}

The results of the previous sections give an algorithm for computing $\PPP(E,E')$ for any rationally 2-isogenous $E$ and $E'$.  Here we explore the range of values it can take.  

\begin{prop} 
With all notation as above, we have $\PPP(E,E') \leq 1/4$.
\end{prop}

\begin{proof}
Suppose $p$ is anomalous for the pair $(E,E')$.  For some $m \ge 2$, we have that either $F \equiv -I \pmod{2^m}$ or $F' \equiv -I \pmod{2^m}$, or both.  Therefore, either $F \equiv -I \pmod{4}$ or $F' \equiv -I \pmod{4}$, which implies that either $-I \in G(4)$ or $-I \in G'(4)$, or both.  We consider two cases.

First suppose that $F \equiv -I \pmod{4}$ for every anomalous prime $p$, or $F' \equiv -I \pmod{4}$ for every anomalous prime $p$.  Switching $E$ and $E'$ as necessary, assume $F \equiv -I \pmod{4}$ for every anomalous prime $p$. In this case, the proportion of anomalous primes is at most 
\[
\lim_{X\rightarrow \infty} \dfrac{\#\{p \leq X: F \equiv -I \pmod{4}\}}{\pi(X)}
\] 
which is precisely ${|G(4)|^{-1}}$ by the Chebotarev Density Theorem. 

We show that $|G(4)|^{-1} \le 1/4$ by showing that $|G(4)| \ge 4$.  The Chebotarev Density Theorem implies that $G(4)$ contains the same number of elements with determinant $1\pmod{4}$ and with determinant $3\pmod{4}$.  Since $\pm I \in G(4)$, two elements with determinant $1$, there must be at least two elements of $G(4)$ with determinant $3$, and so $|G(4)| \ge 4$.

The second case is that for some anomalous prime $p_1,$ $F \not \equiv -I \pmod{4}$; and for some (necessarily distinct) anomalous prime $p_2$, $F' \not \equiv -I \pmod{4}$.  Since $p_1$ (respectively $p_2$) is anomalous, and $F \not \equiv -I \pmod{4}$ (resp.~$F' \not \equiv -I \pmod{4}$), it must be the case that $F' \equiv -I \pmod{4}$ (resp.~$F \equiv -I \pmod{4}$).  Thus, both $G(4)$ and $G'(4)$ contain $-I$, and the pair $(E,E')$ has anomalous primes of defect $(3,2)$ and of defect $(2,3)$.  By Proposition \ref{literal_zero}, these sets have positive Chebotarev densities.

Since both densities are nonzero, we must have $c_4 = c_4' = 1/2$, by Proposition \ref{c2mprop}.  Since each of $G(4)$ and $G'(4)$ contain at least three distinct elements of determinant 1 ($I$, $-I$, $F$, and $I$, $-I$, $F'$, respectively), and since the determinants must be equally distributed modulo 4, it follows that each of $|G(4)|$ and $|G'(4)|$ must be at least $6$.  Note that $|\GL_2(\Z/4\Z)| = 96$.  Since $E$ and $E'$ each have a rational point of order 2, both $G(4)$ and $G'(4)$ have index divisible by $3$ in $\GL_2(\Z/4\Z)$.  This implies that $G(4$) and $G'(4)$ are 2-groups, and so $|G(4)|, |G'(4)| \ge 8$. Therefore, the proportion of primes of defect $(3,2)$ or $(2,3)$ is bounded by
\[ 
\frac{1}{2|G(4)|} + \frac{1}{2|G'(4)|} = \frac{1}{16} + \frac{1}{16} = \frac{1}{8}.
\]
We now must show that the remaining anomalous primes do not contribute more than another 1/8 to this proportion.

If the only anomalous primes for the pair $(E,E')$ have defect $(3,2)$ and $(2,3)$, then we are done by our work above.  If not, then we consider the sizes of the groups $|G(2^m)|$ and $|G'(2^m)|$.  By \cite[Cor.~1.3]{rzb}, the index of $G$ and $G'$ in $\GL_2(\Z_2)$ must divide 96.  

It follows that for each $m \ge 4$, we have both $|G(2^m)|$ and  $|G'(2^m)| \ge 16 \cdot 16^{m-3}$, so in particular both $|G(32)|$ and  $|G'(32)| \ge 16^3$. Therefore, by Theorem \ref{rough_thm}
\[
\PPP(E,E')^{{\rm tail}} \le \frac{8}{15} \left( \frac{1}{|G(32)|}+\frac{1}{|G'(32)|}\right) \le \frac{16}{15\cdot 16^3} = \frac{1}{15\cdot16^2}.
\]
It remains to give an upper bound on the $m=3$ and 4 terms.  

Using the fact (again) that $G$ and $G'$ have index dividing 96 in $\GL_2(\Z_2)$, we have that $|G(8)|, |G'(8)| \geq 16$ and $|G(16)|, |G'(16)| \geq 16^2$.  By Proposition \ref{c2mprop}, we have 
\[
(c_{2^m}, c'_{2^m}) \in \{(0,0), (1/2,0),(0,1/2), (1/2,1/2), (0,1),(1,0)\}.
\]
A maximum theoretical value of $\PPP(E,E')$ is then achieved when all $|G(2^m)|$ and $|G'(2^m)|$  are minimal and weighted by any of the pairs $(c_{2^m}, c'_{2^m}) \in \lbrace (1/2, 1/2), (1,0), (0,1) \rbrace$.  Putting all this together we have 
\begin{align*}
\PPP(E,E') &\leq  \frac{1}{8} + \frac{1}{16} + \frac{1}{16^2} +  \frac{1}{15 \cdot 16^2} \\
&= \frac{23}{120} < \frac{1}{4}.
\end{align*}
\end{proof}

\begin{rmk}
In \cite[Example 5.2.1]{ck} we gave an explicit example where $\PPP(G,G') = 1/4$, showing that this bound can indeed be achieved. Note that in that example, all anomalous primes have defect $(3,2)$.  
\end{rmk}

Turning to small values of $\PPP(E,E')$ the idea is to engineer examples where $c_{2^m} = c_{2^m}'=0$ for small values of $m$, but such that they are non-zero when $G(2^m)$ and $G'(2^m)$ are large.   We present two examples here, due to Jeremy Rouse, of small values of $\PPP(E,E')$ both of which are less than $1/30$.  We thank him for allowing us to include them here.

\begin{exm}
Let $E$ be the elliptic curve \href{https://www.lmfdb.org/EllipticCurve/Q/400/g/1}{400.g1} and $E'$ the elliptic curve \href{https://www.lmfdb.org/EllipticCurve/Q/400/g/2}{400.g2}.  Then $G$ and $G'$ contain $-I$ and each has level 16 with $|G(16)| = |G'(16)| = 256$, so that  
\[
\PPP(E,E')^{\rm tail} = \frac{8}{15} \left(\frac{1}{|G(32)|} + \frac{1}{|G'(32)|} \right) = \frac{8}{15} \left(\frac{1}{4096} + \frac{1}{4096} \right) = \frac{1}{3840}.
\]
Explicit computations with the division polynomials in \textsf{Magma} show that the 4-torsion and 8-torsion fields of $E$ and $E'$ are isomorphic, whence 
$c_4 = c_4' = c_8 = c_8' = 0$. 
However, for the 16-torsion we are in Case 3(b) and 4(b) of Proposition \ref{c2mprop} so that $c_{16} = c_{16}' = 1/2$ and so $\PPP(E,E')^{\rm head} = (1/2)/256 + (1/2)/256 = 1/256$.  Therefore, 
\[
\PPP(E,E') = \PPP(E,E')^{\rm head}  + \PPP(E,E')^{\rm tail} = \frac{1}{256} + \frac{1}{3840} = \frac{1}{240}.
\]

Using \textsf{Magma} we checked these numbers experimentally.  We sampled the first 100,000 good primes for $E$ and $E'$ and found
\begin{center}
\begin{tabular}{rl}
192 &\text{primes of defect (4,5)} \\
185 &\text{primes of defect (5,4)} \\
15 &\text{primes of defect (5,6)} \\
15 &\text{primes of defect (6,5)} \\
1 &\text{prime of defect (6,7)} \\
\hline
\textbf{408} & anomalous primes
\end{tabular}
\end{center}
Note that based on our prediction above, we would expect a total of $100,000/240 \approx 417$ anomalous primes, with
\[
\frac{1}{2} \cdot \frac{1}{256} \cdot 100,000 \approx 195
\]
primes of defect (4,5) and (5,4), each and
\[
\frac{1}{2} \cdot \frac{1}{4096} \cdot 100,000 \approx 12
\]
primes of defect (5,6) and (6,5), each.
\end{exm}

\begin{exm}
Let $E$ be the elliptic curve \href{https://www.lmfdb.org/EllipticCurve/Q/4225/m/1}{4225.m1} and $E'$ the elliptic curve \href{https://www.lmfdb.org/EllipticCurve/Q/4225/m/2}{4225.m2}.  Then $G$ and $G'$ contain $-I$ and each has level 16 with $|G(16)| = |G'(16)| = 1024$, so that  
\[
\PPP(E,E')^{\rm tail} = \frac{8}{15} \left(\frac{1}{|G(32)|} + \frac{1}{|G'(32)|} \right) = \frac{8}{15} \left(\frac{1}{16384} + \frac{1}{16384} \right) = \frac{1}{15360}.
\]
The rest of the situation is identical to the previous example: the 4-torsion and 8-torsion fields of $E$ and $E'$ are isomorphic, whence $c_4 = c_4' = c_8 = c_8' = 0$. We compute that $c_{16} = c_{16}' = 1/2$, and so $\PPP(E,E')^{\rm head} = 1/1024$.  Altogether, we have $\PPP(E,E') = 1/960$.

Turning to \textsf{Magma}, we sample the first 100,000 good primes and find
\begin{center}
\begin{tabular}{rl}
50 &\text{primes of defect (4,5)} \\
42 &\text{primes of defect (5,4)} \\
2 &\text{primes of defect (5,6)} \\
3 &\text{primes of defect (6,5)} \\
1 &\text{prime of defect (6,7)} \\
\hline
\textbf{98} & anomalous primes
\end{tabular}
\end{center}
Based on our prediction above, we would expect a total of $100,000/960 \approx 104$ anomalous primes, with
\[
\frac{1}{2} \cdot \frac{1}{1024} \cdot 100,000 \approx 49
\]
primes of defect (4,5) and (5,4), each, and
\[
\frac{1}{2} \cdot \frac{1}{16384} \cdot 100,000 \approx 3
\]
primes of defect (5,6) and (6,5), each.
\end{exm}

\section{Elliptic Curves with CM} \label{cmsection}

When $E$ and $E'$ have complex multiplication over $\Q$, the situation is much simpler, owing to the fact that the CM structure of $E$ and $E'$ over $\Q$ is preserved modulo $p$.  The goal of this section is to prove Theorem \ref{cm_thm}, showing that when $E$ and $E'$ have CM there are only two possibilities for $\PPP(E,E')$; either it is 0 or it is 1/12.

Let $E$ be a curve with CM by an order $\mathcal{O}$ in a quadratic imaginary field and let $p$ be a prime of good reduction for $E$. By \cite[Ch.~13, Thm.~12]{lang}, $E$ has ordinary reduction if and only if $p$ splits in $K$ and, outside of the finitely many primes dividing the conductor, we have $\End(E \pmod{p}) \cong \mathcal{O}$. Thus (writing $\mathcal{O}'$ for the endomorphism ring of $E'$),  whenever the reduction is ordinary, the reduction of $E \to E'$ will (outside of finitely many exceptions) be horizontal (if $\mathcal{O} \cong \mathcal{O}'$), ascending (if $\mathcal{O} \subsetneq \mathcal{O}'$), or descending (if $\mathcal{O} \supsetneq \mathcal{O}'$). This is in contrast to the non-CM case, especially when the groups $G$ and $G'$ are maximal and the rational isogeny $E \to E'$ appears to reduce modulo $p$ to an ascending/descending/horizontal isogeny ``at random'' (see \cite[\S 6]{ck} for a precise statement).  

The key starting point for our result on CM elliptic curves is Proposition \ref{descending_prop} below.  For completeness, we remind the reader of several exceptional cases which do not affect our study of anomalous primes:
\begin{enumerate}
\item If $p$ is a good prime for which the reduction of $E$ and $E'$ is supersingular, then $p$ cannot be anomalous, by \cite[Thm.~2.3.1]{ck}.
\item The endomorphism rings will not be preserved under reduction modulo $p$ if $p$ divides the conductor of $\mathcal{O}$ or $\mathcal{O}'$.  Indeed, by \cite[Ch.~13, Thm.~12]{lang}, if $\mathcal{O}$ has conductor $f = p^r f'$ where $p$ does not divide $f'$, then the endomorphism ring of the reduction, say $\mathcal{O}_p$, has conductor $f'$.  But since $E$ and $E'$ are defined over $\Q$, they each have rational $j$-invariant.  Thus $\mathcal{O}$ and $\mathcal{O}'$ have class number 1,  which means the conductors of $\mathcal{O}$ and $\mathcal{O}'$ are in $\{1,2,3\}$.

Let  $\mathcal{O}_p$ and $\mathcal{O}'_p$ denote the endomorphism rings of the reduction at $p$.  If $p$ is anomalous, then $\mathcal{O}_p$ and $\mathcal{O}'_p$ must have relative index 2. 
That means either $\mathcal{O}$ and $\mathcal{O}'$ have relative index 2 and $p$ divides the conductor or each, or $\mathcal{O}$ and $\mathcal{O}'$ have relative index larger than and divisible by 2.  Neither one of these situations is possible in our setup.  We conclude that if $p$ divides the conductor of $\mathcal{O}$ or $\mathcal{O}'$, $p$ is necessarily not anomalous. 
\end{enumerate}
The upshot of these two points is that we can safely ignore the (infinitely many) primes of supersingular reduction and the (finitely many) primes dividing the conductors of $\mathcal{O}$ and $\mathcal{O}'$.

\begin{prop} \label{descending_prop}
Let $E$ and $E'$ be rationally 2-isogenous elliptic curves defined over $\Q$ with CM by the orders $\mathcal{O}$ and $\mathcal{O}'$, respectively. Then exactly one of the following is true:
\begin{itemize}
\item 
if $p$ is an anomalous prime, then $p$ has defect $(m+1, m)$ for some $m \ge 2$, or
\item 
if $p$ is an anomalous prime, then $p$ has defect $(m,m+1)$ for some $m \ge 2$, or
\item there are no anomalous primes.
\end{itemize}
\end{prop}

\begin{proof}
Let $p>3$ be a good prime for which the reduction is ordinary.   Then modulo $p$, the elliptic curves have $\End(E) = \mathcal{O}$ and $\End(E') = \mathcal{O}'$ as well, by \cite[Ch.~13, Thm.~12]{lang}.   If $\mathcal{O} = \mathcal{O}'$, the group structures $E(\F_{p^k})$ and $E'(\F_{p^k})$ are isomorphic for all $k\geq 1$ and hence $p$ is not anomalous.  Otherwise, we have $[\mathcal{O}:\mathcal{O}'] = 2^{\pm 1}$ . Thus, the isogeny $E \to E'$ is either descending for all good primes or is ascending for all good primes.  Therefore, if the set of anomalous primes is non-empty, then there are exactly two possibilities: either every anomalous prime is of defect $(m+1,m)$ for some $m \ge 2$, or every anomalous prime is of the form $(m,m+1)$ for some $m \ge 2$. 
\end{proof}

We now turn to two classifications of CM elliptic curves over $\Q$ that, taken together, reduce the number of isomorphism types of $G$ and $G'$ attached to $E$ and $E'$ to a total of four pairs that have a nonzero proportion of anomalous primes.  In \cite{alvaro-chiloyan}, Chiloyan and Lozano-Robledo classify all possible isogeny-torsion graphs of elliptic curves over $\Q$ and also indicate which of these graphs occur for curves with CM. In \cite{alvaro}, Lozano-Robledo determines all possible 2-adic images attached to CM elliptic curves over $\Q$. Thus, \cite{alvaro} tells us which images occur, and \cite{alvaro-chiloyan} tells us how they can be connected by 2-isogenies.

Of the 52 possible isogeny-torsion graphs of elliptic curves over $\Q$, only the graphs $L_2(2)$, $T_4$ (including three sub-cases), $R_4(6)$ (including two sub-cases), and $R_4(14)$ can occur for elliptic curves over $\Q$ with CM that also have a rational 2-torsion point.  When $E$ and $E'$ are rationally 2-isogenous elliptic curves over $\Q$ with CM, there are seven possible isomorphism types of isogeny-torsion graphs on which they can represent adjacent vertices. These graphs are listed in Figure \ref{isotor}, with unlabeled edges representing isogenies of degree 2; otherwise we label edges by the degree of the isogeny.

\begin{figure}[htbp]
\begin{center}
\begin{tabular}{|c|c|}
\hline
Label & Graph\\
\hline
$L_2(2)$ & \xymatrix{[2] \ar@{-}[r] & [2]}\\ 
\hline
$T_4^1$ & $\xymatrix{
&[4] \ar@{-}[d] \\
[4]&[2,2] \ar@{-}[l] \ar@{-}[r] & [2] 
}$ \\
\hline
$T_4^2$ & $\xymatrix{
&[2] \ar@{-}[d] \\
[4]&[2,2] \ar@{-}[l] \ar@{-}[r] & [2] 
}$ \\
\hline
$T_4^3$ & $\xymatrix{
&[2] \ar@{-}[d] \\
[2]&[2,2] \ar@{-}[l] \ar@{-}[r] & [2] 
}$  \\
\hline
$R_4(6)$ &$\xymatrix{
[2] \ar@{-}[r] & [2] \ar@{-}^3[d] \\
[6] \ar@{-}^3[u] \ar@{-}[r] & [6]}$  \\
\hline
$R_4(6)$ &$\xymatrix{
[2] \ar@{-}[r] & [2] \ar@{-}^3[d] \\
[2] \ar@{-}^3[u] \ar@{-}[r] & [2]}$  \\
\hline
$R_4(14)$ &$\xymatrix{
[2] \ar@{-}[r] & [2] \ar@{-}^7[d] \\
[2] \ar@{-}^7[u] \ar@{-}[r] & [2]}$\\
\hline
\end{tabular}
\end{center}
\caption{Isogeny-Torsion Graphs of CM curves with a point of order 2} \label{isotor}
\end{figure}

For each isogeny-torsion graph in Figure \ref{isotor}, there may be multiple possible isomorphism types for the images of the 2-adic representations of the elliptic curves corresponding to the vertices of the graph. 
We will now systematically go through all of the possible cases of 2-adic images as classified in \cite{alvaro}.  As a guide to the classification, we note 
\begin{itemize}
\item images for curves with $j$-invariant 1728 are classified in \cite[Thm.~9.7]{alvaro}, and
\item images for curves with $j$-invariant 0 are classified in \cite[Thm.~9.10]{alvaro}, and
\item all other images are classified in \cite[Thm.~9.3]{alvaro}.
\end{itemize}

\subsection{Step 1 -- Classify images for which there are no anomalous primes}

Using both \cite{alvaro} and \cite{alvaro-chiloyan} as organizational guides, we will first eliminate from consideration all rationally 2-isogenous elliptic curves over $\Q$ with CM for which there are no anomalous primes.  We start with a general lemma before working through specific cases.  Throughout the rest of this section we always assume that $E$ has CM, which means that any curves isogenous to it also have CM, even if we do not explicitly recall this assumption every time we mention~$E$. 

\begin{lem} \label{equal_j}
If $E$ and $E'$ are rationally 2-isogenous elliptic curves with equal $j$-invariants, then there are no anomalous primes for $(E,E')$.
\end{lem}

\begin{proof}
Since $j(E) = j(E')$, the curves $E$ and $E'$ are isomorphic over $\overline{\Q}$.  This implies that for every good prime $p$, their reductions modulo $p$ are isomorphic over $\overline{\F}_p$, and, in fact, they are isomorphic over some finite extension $\F_{p^k}$ of $\F_p$.  Thus $\End_{\F_{p^k}}(E) = \End_{\F_{p^k}}(E')$.  But since $E$ and $E'$ have ordinary reduction at $p$, all endomorphisms are defined over $\F_p$, hence $\mathcal{O} = \mathcal{O'}$, proving the lemma.
\end{proof}

\noindent Now we turn to specific cases.

\subsubsection{$j$-invariant 1728} These elliptic curves can be taken to have Weierstrass form 
\[
E: y^2 = x^3 + a_4x.
\]
Let $A = (x_A,0)$ be a 2-torsion point of $E$, whether defined over $\Q$ or not.  Then V\'elu's formulas \cite{velu} imply that $E' = E/\langle A \rangle$ has Weierstrass form
\[
E' : y^2 = x^3 +\left( -15x_A^2 - 4a_4\right)x + 14a_4x_A.
\]
Since the point $(0,0)$ is $\Q$-rational, $E$ always admits a rational 2-isogeny.  When neither $a_4$ nor $-a_4$ is a rational square, the curves $E$ and $E'$ belong to an isogeny-torsion graph of type $L_2(2)$.  Otherwise, one of $\pm a_4$ is a rational square, and $E$ and $E'$ make two vertices of an isogeny-torsion graph of type $T_4^3$ \cite[Thm.~9.7]{alvaro}.

\begin{lem}
If neither $a_4$ nor $-a_4$ is a rational square, then $E(\Q)[2^\infty] \simeq E'(\Q)[2^\infty] \simeq \Z/2\Z$ and there are no anomalous primes for this pair.
\end{lem}

\begin{proof}
In this case the curves $E$ and $E'$ lie on an isogeny-torsion graph of type $L_2(2)$ and each has $j$-invariant 1728. By Lemma \ref{equal_j} there are no anomalous primes.
\end{proof}

If $a_4 = -t^2$, then the curves
\begin{align*}
E:~& y^2 = x^3 -t^2x \\
E':~&y^2 = x^3 + 4t^2x \\
E'':~&y^2 = x^3 - 11t^2x -14t^3 \\
E''':~&y^2 = x^3 - 11t^2x + 14t^3,
\end{align*}
form an isogeny-torsion graph of type $T_4^3$. It is a computation with 2- and 4-torsion polynomials to check that
\[
E'(\Q)_\tors[2^\infty] \simeq E''(\Q)_\tors[2^\infty] \simeq E'''(\Q)_\tors[2^\infty] \simeq \Z/2\Z,
\] 
while $E(\Q)_\tors[2^\infty] \simeq \Z/2\Z \times \Z/2\Z$.  We have $j(E) = j(E') = 1728$ and $j(E'') = j(E''') = 287496$.   Similarly, if $a_4 = t^2$, then the curves
\begin{align*}
E:~& y^2 = x^3 +t^2x \\
E':~&y^2 = x^3 - 4t^2x \\
E'':~&y^2 = x^3 -44t^2 -112t^3 \\
E''':~&y^2 = x^3  -44t^2 +112t^3,
\end{align*}
form an isogeny-torsion graph of type $T_4^3$.  In this case, we compute that 
\[
E'(\Q)[2^\infty] \simeq \Z/2\Z \times \Z/2\Z,\ \ \ \text{ and } \ \ \ E(\Q)[2^\infty] \simeq E''(\Q)[2^\infty]  \simeq E'''(\Q)[2^\infty] \simeq \Z/2\Z,
\]
and $j(E) = j(E') = 1728$, $j(E'') = j(E''') = 287496$. By Lemma \ref{equal_j} we need only consider the 2-isogenies linking the curves of distinct $j$-invariants.  We put this off until Subsection \ref{1/12} below.

\subsubsection{$j$-invariant 0} 
An elliptic curve $E$ of $j$-invariant $0$ has a short Weierstrass form
\[
E: y^2 = x^3 + a_6.
\]
Importing the notation of \cite[Thm.~9.10]{alvaro}, there are two possibilities for 2-adic image attached to such a curve:
\[
[\mathcal{C}_{-1,1}(2^\infty):G_{E,K,2^\infty}]  \in \lbrace 1, 3 \rbrace.
\]
However, elliptic curves for which $[\mathcal{C}_{-1,1}(2^\infty):G_{E,K,2^\infty}]  =1$ do not have a $\Q$-rational 2-torsion point and so do not admit a rational 2-isogeny.  

We have $[\mathcal{C}_{-1,1}(2^\infty):G_{E,K,2^\infty}]  =3$ if and only if $a_6$ is a rational cube and $E(\Q)[2^\infty] \simeq \Z/2\Z$.  Thus $E$ admits exactly one rational 2-isogeny $E \to E'$, where
\begin{align*}
E:~& y^2 = x^3 + t^3 \\
E':~&y^2 = x^3 -15t^2x+22t^3,
\end{align*}
by V\'elu's formulas; it is routine to check that $E(\Q)[2^\infty] \simeq E'(\Q)[2^\infty] \simeq \Z/2\Z$, $j(E) = 0$, and $j(E') = 54000$.  Consider also the twists
\begin{align*}
E'':~&y^2 = x^3 - 27t^3 \\
E''':~&y^2 = x^3 -135t^2 -594t^3;
\end{align*}
we similarly have $E''(\Q)[2^\infty] \simeq E'''(\Q)[2^\infty] \simeq \Z/2\Z$, $j(E'') = 0$, and $j(E''') = 54000$.  Then the curves $E,E',E'',E'''$ fit into an isogeny-torsion graph of type $R_4(6)$
\[
\xymatrix{
E \ar@{-}^3[r] \ar@{-}_2[d]&  E'' \ar@{-}^2[d] \\
E'  \ar@{-}^3[r] & E'''
}
\]
where either all curves have torsion subgroup $\Z/2\Z$, or exactly one of $\lbrace E,E'' \rbrace$ and exactly one of $\lbrace E', E''' \rbrace$ have torsion subgroup $\Z/6\Z$ while the others have torsion subgroup $\Z/2\Z$.  We further investigate the 2-isogenous curves in Subsection \ref{1/12} below. 

\subsection{Isogeny-Torsion graphs with a 4-torsion point} 
By Lemma \ref{4lem}, any pair of 2-isogenous elliptic curves for which at least one of them has a rational 4-torsion point cannot have any anomalous primes.  This affects elliptic curves that belong to the isogeny-torsion graphs $T_4^1$ and $T_4^2$. 

\begin{prop} \label{T_4^1-graph}
Let $E$ and $E'$ be rationally 2-isogenous elliptic curves belonging to an isogeny-torsion graph of type $T_4^1$.  Then they have no anomalous primes.
\end{prop}

\begin{proof}
It suffices to check the subgraph $\xymatrix{[2,2] \ar@{-}[r] & [2]}$.  Without loss of generality suppose $E(\Q)[2] \simeq \Z/2\Z \times \Z/2\Z$ and $E'(\Q)[2] \simeq \Z/2\Z$.  By the main classification result in \cite{alvaro}, using the notation from that paper, the full 2-adic images are
\begin{align*}
G &\simeq \left< \begin{pmatrix} -1 & 0 \\ 0 & 1 \end{pmatrix}, \begin{pmatrix} 5 & 0 \\ 0 & 5\end{pmatrix}, \begin{pmatrix} -1 & -2  \\ 2 & -1 \end{pmatrix} \right> =  \langle c_{-1}, G_{4,b} \rangle, \text{ and} \\
G'&\simeq \left< \begin{pmatrix} -1 & 0 \\ 0 & 1 \end{pmatrix}, \begin{pmatrix} 5 & 0 \\ 0 & 5\end{pmatrix}, \begin{pmatrix} -1 & -1  \\ -4 & -1 \end{pmatrix} \right>.
\end{align*} 
It is a computation to check that neither $G$ nor $G'$ contains $-I \pmod{4}$. Hence there exist no primes of defect $(m+1,m)$ or $(m,m+1)$ for all $m \geq 2$, and hence no anomalous primes.
\end{proof}

\begin{prop}
Let $E$ and $E'$ be rationally 2-isogenous elliptic curves belonging to an isogeny-torsion graph of type $T_4^2$.  Then they have no anomalous primes.
\end{prop}

\begin{proof}
By the classification of Galois images associated to CM elliptic curves from \cite{alvaro}, 
type $T_4^2$ isogeny-torsion graphs occur under the following conditions.  We have elliptic curves $E,E',E'',E'''$ with center vertex $E$; let $j,j',j'',j'''$ denote their $j$-invariants, respectively.  In this setup we have $j = j'$ and $j'' = j'''$, where $E'(\Q)[2^\infty] \simeq E''(\Q)[2^\infty] \simeq \Z/2\Z$ and $E'''(\Q)[2^\infty] \simeq \Z/4\Z$.  

By Lemma \ref{4lem} the rationally 2-isogenous pair $(E,E''')$ cannot have any anomalous primes.  Since $j=j'$, Lemma \ref{equal_j}  implies that $(E,E')$ cannot have any anomalous primes.  
It remains to check the pair $(E,E'')$.  The proof is similar to the proof of Proposition \ref{T_4^1-graph}.  According to the classification of \cite{alvaro}:
\begin{align*}
G &\simeq \left< \begin{pmatrix} -1 & 0 \\ 0 & 1 \end{pmatrix}, \begin{pmatrix} 5 & 0 \\ 0 & 5\end{pmatrix}, \begin{pmatrix} 1 & 2  \\ -2 & 1 \end{pmatrix} \right> =  \langle c_{-1}, G_{4,a} \rangle, \text{ and} \\
G''&\simeq \left< \begin{pmatrix} -1 & 0 \\ 0 & 1 \end{pmatrix}, \begin{pmatrix} 5 & 0 \\ 0 & 5\end{pmatrix}, \begin{pmatrix} 1 & 1  \\ -4 & 1 \end{pmatrix} \right>.
\end{align*} 
It is a computation to check that neither $G$ nor $G'$ contains $-I\pmod{4}$.  We conclude that there are no anomalous primes in this case.
\end{proof}

\subsubsection{Isogeny-Torsion Graphs of type $L_2(2)$} Now we consider the two final cases for which there are no anomalous primes.

\begin{prop} \label{L_2(2)}
Suppose $E$ and $E'$ are rationally 2-isogenous elliptic curves with CM and $j$-invariants not equal to $0$ or $1728$ that belong to an isogeny-torsion graph of type $L_2(2)$.  Then there are no anomalous primes for this pair.
\end{prop}

\begin{proof}
There are only two cases to consider.  According to the classification of \cite[Thm.~9.3]{alvaro}: either 
\begin{align*}
G &\simeq \left< \begin{pmatrix} 1 & 0 \\ 0 & -1 \end{pmatrix}, \begin{pmatrix} 3 & 0 \\ 0 & 3\end{pmatrix}, \begin{pmatrix} -1 & -1  \\ 4 & -1 \end{pmatrix} \right>, \text{ and} \\
G''&\simeq \left< \begin{pmatrix} -1 & 0 \\ 0 & 1 \end{pmatrix}, \begin{pmatrix} 3 & 0 \\ 0 & 3\end{pmatrix}, \begin{pmatrix} -1 & -1  \\ 4 & -1 \end{pmatrix} \right>,
\end{align*} 
or
\begin{align*}
G &\simeq \left< \begin{pmatrix} 1 & 0 \\ 0 & -1 \end{pmatrix}, \begin{pmatrix} 3 & 0 \\ 0 & 3\end{pmatrix}, \begin{pmatrix} 1 & 1  \\ -4 & 1 \end{pmatrix} \right>, \text{ and} \\
G''&\simeq\left< \begin{pmatrix} -1 & 0 \\ 0 & 1 \end{pmatrix}, \begin{pmatrix} 3 & 0 \\ 0 & 3\end{pmatrix}, \begin{pmatrix} 1 & 1  \\ -4 & 1 \end{pmatrix} \right>.
\end{align*} 
In all four possibilities the corresponding $j$-invariants are equal to $2^65^3$ by \cite[Example 9.4 (2)]{alvaro}. By Lemma \ref{equal_j} there are no anomalous primes. 
\end{proof}

\subsection{Anomalous Primes in the CM Case} \label{1/12}

We are now left with a handful of cases, all of which (as we will show in this section) do have anomalous primes.  In each case we will show that the proportion of anomalous primes is 1/12.  To start, we list the remaining cases of pairs of rationally 2-isogenous elliptic curves by the isomorphism types of their 2-adic images. We mark which curves have $j$-invariant 0 or 1728 so that the reader can cross-check these cases with those from earlier in this section.
\begin{center}
\begin{tabular}{|c|c|c|}
\hline
Isogeny-Torsion Graph & 2-adic images& Example \\
\hline
$R_4(14)$ & $\xymatrix{N_{\delta,\phi}(2^\infty) \ar@{-}[r] & N_{\delta,\phi}(2^\infty)}$ &  \href{https://www.lmfdb.org/EllipticCurve/Q/49/a/}{49.a}\\

\hline
$R_4(6)$ &  $\xymatrix{N_{\delta,\phi}(2^\infty)  \ar@{-}[r] & \langle C_{-1,1}(2^\infty)^3,c_{-1}' \rangle(j=0)} $ & \href{https://www.lmfdb.org/EllipticCurve/Q/36/a/}{36.a} \\
\hline
$R_4(6)$ &  $\xymatrix{N_{\delta,\phi}(2^\infty)   \ar@{-}[r] & \langle C_{-1,1}(2^\infty)^3,c_{1}'  \rangle (j=0)}$ & \href{https://www.lmfdb.org/EllipticCurve/Q/144/a/}{144.a}\\
\hline
$T_4^3$ & $\xymatrix{ &  N_{\delta,\phi}(2^\infty) \ar@{-}[d] & \\ 
N_{\delta,\phi}(2^\infty)  \ar@{-}[r]& \langle G_{2,a},c_{-1} \rangle \atop{j=1728} &  \ar@{-}[l]\langle G_{2,a},c_{-1}'\rangle \atop{j=1728} & }$  & \href{https://www.lmfdb.org/EllipticCurve/Q/288/d/}{288.d}\\
\hline
\end{tabular}
\end{center}

\begin{rmk}
In the two isogeny-torsion graphs marked $R_4(6)$ in the table above, the elliptic curves of $j$-invariant 0 have $a_4=0$ and $a_6 = t^3$ for some $t \in \Q^\times$ \cite[Thm.~9.3]{alvaro}.  If $(-t,0)$ denotes the 2-torsion point of $E$ generating the kernel of the unique 2-isogeny $E \to E'$, then we have
\begin{align*}
E:~&y^2 = x^3 + t^3 \\
E':~&y = x^3 - 15t^2x + 22t^3,
\end{align*}
and $j(E') = 54000$. 
\end{rmk}

Now we turn to the proof of Theorem \ref{cm_thm}:
if $E$ and $E'$ are rationally 2-isogenous elliptic curves with CM then either $\PPP(E,E') = 0$ or $\PPP(E,E') = 1/12$.  According to the table above, the only remaining cases for which $\PPP(E,E')$ is potentially nonzero are:
\begin{align*}
&\xymatrix{N_{\delta,\phi}(2^\infty) \ar@{-}[r] & N_{\delta,\phi}(2^\infty),} \\
&\xymatrix{N_{\delta,\phi}(2^\infty)  \ar@{-}[r] & \langle C_{-1,1}(2^\infty)^3,c_{-1}'\rangle,} \\
&\xymatrix{N_{\delta,\phi}(2^\infty)  \ar@{-}[r] & \langle C_{-1,1}(2^\infty)^3,c_{1}'\rangle, \text{ and}} \\
&\xymatrix{N_{\delta,\phi}(2^\infty)  \ar@{-}[r]& \langle G_{2,a},c_{-1} \rangle.}
\end{align*}
The following is easily verifiable based on the explicit descriptions in \cite[\S 9]{alvaro}, and we omit the proof.

\begin{lem} \label{CMkernel}
Let $G \in \lbrace N_{\delta,\phi}(2^\infty), \langle C_{-1,1}(2^\infty)^3,c_{-1}'\rangle, \langle C_{-1,1}(2^\infty)^3,c_{1}'\rangle, \langle G_{2,a},c_{-1} \rangle \rbrace$.  Then for any $m \ge 2,\ |G(2^m)| = 2 \cdot 4^{m-1}$. 
\end{lem}

In each of these four cases, our strategy for showing that $\PPP(E,E') = 1/12$ will be the same, and is influenced by the proof of \cite[Thm.~5.2.4]{ck}.  

\subsection{Interlude on Isogeny Volcanoes} \label{volcano_section}  A main ingredient in the proof of Theorem \ref{cm_thm} is the structure of the isogeny volcano corresponding to an elliptic curve over a finite field.   In this section we will collect all of the relevant formulas and ideas used in the proof in a concise manner.   We refer the reader to \cite[\S4]{ck} for a more thorough introduction or \cite{sutherland_volcano} for an extensive treatment.

The 2-isogeny volcano $V_q$ corresponding to an elliptic curve $E$ over a finite field $\Fq$ is a graph whose vertices correspond to isomorphism classes of elliptic curves over $\Fq$ and whose edges are 2-isogenies.  The vertices are partitioned into subsets called \textbf{levels}, and elliptic curves at the same level have isomorphic endomorphism rings and isomorphic 2-Sylow subgroups; the levels are indexed from $0$ (the floor) to $h(V_q)$ (the crater); the quantity $h(V_q)$ is called the \textbf{height} of the volcano, which we next describe how to compute.

Let $t$ be the trace of Frobenius on $E$ over $\F_q$ and, for a non-zero integer $n$, write $\sqf(n)$ for the squarefree part of $n$.  Set $D = \sqf(t^2-4q)$ so that $\End(E) \otimes \Q \simeq K \ddef \Q(\sqrt{D})$.  We write $h(V_q)$ for the height of $V_q$; by \cite[Thm.~25]{kohel} we have 
\[
h(V_q) = \frac{1}{2} v_2 \left( \frac{t^2 - 4q}{\disc \mathcal{O}_0} \right) = \frac{1}{2} v_2 \left(\frac{t^2 - 4q}{\disc \OK} \right),
\]
where $\mathcal{O}_0$ is the endomorphism ring of an elliptic curve lying on the crater of $V_q$.  By \cite[Thm.~7(5)]{kohel}, $[\OK:\mathcal{O}_0]$ is odd.

Because the elliptic curves belonging to $V_q$ are isogenous, all 2-Sylow subgroups have the same size.   The 2-Sylow subgroups at level 0 are cyclic.  A 2-isogeny volcano is typically decorated with the 2-Sylow subgroup at each level.

Now set $q=p$ and suppose $p$ is an anomalous prime for the pair $(E,E')$.  Then the floor of $V_p$ consists of vertices decorated with the group $\Z/4\Z$, while every other level is decorated with $\Z/2 \Z \times \Z/2\Z$.  Upon base change to $\F_{p^2}$, if $|E(\F_{p^2})[2^\infty]| = 2^v$, then (depending on the height of $V_{p^2}$), the group structures ``balance'' as the level increases:
\[
\underbrace{\Z/2^v \Z}_{\text{level 0}}, \ \underbrace{\Z/2 \times \Z/2^{v-1}\Z}_{\text{level 1}}, \ \underbrace{\Z/2^2 \times \Z/2^{v-2}\Z}_{\text{level 2}}, \ \dots.
\]
For example, we might see anomalous prime behavior in the levels of a volcano as follows:
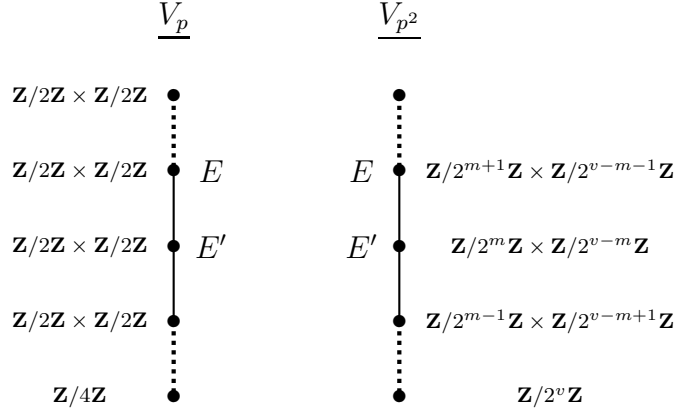
\begin{center}
\begin{tikzpicture}[scale=1.0,sizefont/.style={scale = 2}]
\draw[ultra thick] (1,5) node {${\bullet}$};
\draw[ultra thick, dotted] (1,5) -- (1,4);
\draw[ultra thick] (1,4) node {${\bullet}$};
\draw[thick] (1,3) -- (1,4);
\draw[ultra thick] (1,3) node {${\bullet}$};
\draw[ultra thick] (1,2) node {${\bullet}$};
\draw[ultra thick, dotted] (1,2) -- (1,1);
\draw[ultra thick] (1,1) node {${\bullet}$};
\draw[thick] (1,3) -- (1,2);
\draw(6,5) node {};
\draw(6,4) node {\tiny{$\Z/2^{m+1}\Z \times \Z/2^{v-m-1}\Z$}};
\draw(6,3) node {\tiny{$\Z/2^m\Z \times \Z/2^{v-m}\Z$}};
\draw(6,2) node {\tiny{$\Z/2^{m-1}\Z \times \Z/2^{v-m+1}\Z$}};
\draw(6,1) node {\tiny{$\Z/2^v\Z$}};
\draw(3.5,4) node {$E$};
\draw(3.5,3) node {$E'$};
\draw[ultra thick] (4,5) node {${\bullet}$};
\draw[ultra thick, dotted] (4,5) -- (4,4);
\draw[ultra thick] (4,4) node {${\bullet}$};
\draw[thick] (4,3) -- (4,4);
\draw[ultra thick] (4,3) node {${\bullet}$};
\draw[ultra thick] (4,2) node {${\bullet}$};
\draw[ultra thick, dotted] (4,2) -- (4,1);
\draw[ultra thick] (4,1) node {${\bullet}$};
\draw[thick] (4,3) -- (4,2);
\draw(-0.25,5) node {\tiny{$\Z/2\Z \times \Z/2\Z$}};
\draw(-0.25,4) node {\tiny{$\Z/2\Z \times \Z/2\Z$}};
\draw(-0.25,3) node {\tiny{$\Z/2\Z \times \Z/2\Z$}};
\draw(-0.25,2) node {\tiny{$\Z/2\Z \times \Z/2\Z$}};
\draw(-0.25,1) node {\tiny{$\Z/4\Z$}};
\draw(1.5,4) node {$E$};
\draw(1.5,3) node {$E'$};
\draw(4,6) node {\underline{$V_{p^2}$}};
\draw(1,6) node {\underline{$V_{p}$}};
\end{tikzpicture}
\end{center}
If $p$ is an anomalous primes for the pair $(E,E')$, then we know that $t \equiv 2 \pmod{2}$ and \cite[Lem.~4.0.2]{ck} implies that $h(V_{p^2}) = h(V_{p}) +1$. All of this structure gives us an avenue for investigating anomalous primes.  

We prove Theorem \ref{cm_thm} by showing that in every remaining possibility for CM curves over $\Q$, there exist anomalous primes of defect $(m+1,m)$ for all $m \geq 2$.  To do this, we work explicitly with the matrix groups enumerated in Lemma \ref{CMkernel}; specifically, we fix $m \geq 2$ and exhibit a matrix $F \equiv -I \pmod{2^m}$ belonging to $G(2^m)$ such that if $F$ is in the conjugacy class of Frobenius for $E$, then $E$ lies at the correct level to ensure that $E(\F_{p^2})[2^\infty] \not\simeq E'(\F_{p^2})[2^\infty]$ much like in the figure above.  Note that if $h(V_{p^2})$ is large compared to $v_2(|E(\F_{p^2})|)$, then the group structures at the higher levels of $V_{p^2}$ may become isomorphic.  It is therefore imperative that we understand how $v_2(|E(\F_{p^2})|)$, the level of $E$, and the height of $V_p$  are related. 

This argument requires some additional ingredients from the theory of isogeny volcanoes.  If $p$ is anomalous of defect $(m+1,m)$ and $F \in \GL_2(\Z_2)$ is in the conjugacy class of Frobenius, then 
\[
F = -I + 2^m \left(\begin{array}{cc} x & y \\ z & w
\end{array}\right)
\]
where $x,y,z,w$ are not all $0\pmod{2}$.  The discriminant  $\disc \mathcal{O}_0 \pmod{8}$ is determined by $\sqf((x-w)^2 + 4yz) \pmod{8}$ and the height of $V_p$ is $H$ if and only if 
\begin{align} \label{height_matrix}
v_2(t^2 -4p) = 2m + v_2((x-w)^2 + 4yz) =  \begin{cases} 2H & \text{if }  D\equiv 1 \pmod{4} \\
2H+2 & \text{if } D\equiv 3 \pmod{4} \\
2H+3 & \text{if } D\equiv 2 \pmod{4}
\end{cases}.
\end{align}

\medskip

We will use the following observation multiple times in the proof of Theorem \ref{cm_thm}, so for convenience we present it here as a proposition.

\begin{prop} \label{key_height_prop}
Let $\varphi: E \to E'$ be a $\Q$-rational 2-isogeny.  Suppose $p$ is a good prime for $E$ and $E'$ such that $\varphi$ reduces modulo $p$ to a descending isogeny.  Fix a positive integer $m \ge 2$ and suppose $F \equiv -I \pmod{2^m}$.  If $E$ is at level $m$ on $V_p$, then $p$ is anomalous of defect $(m+1,m)$.
\end{prop}

\begin{proof}
Because $F \equiv -I \pmod{2^m}$, we have $E(\F_p)[2^\infty] \simeq \Z/2\Z \times \Z/2\Z$ and $p \equiv 1 \pmod{2^m}$.  Therefore, $t \equiv 2\pmod{4}$ and so $h(V_{p^2}) = h(V_{p}) +1$.  In particular, if $E$ is at level $m$ on $V_p$, then $E$ is at level $m+1$ on $V_{p^2}$.

Because $F \equiv -I \pmod{2^m}$, it follows that $E(\F_{p^2})[2^\infty] \supseteq \Z/2^{m+1}\Z \times \Z/2^{m+1}\Z$.  Suppose $v_2(|E(\F_{p^2})|) = v$.  Since $E \to E'$ is descending, we have enough information to reconstruct the group structures of the first $m+1$ levels of $V_{p^2}$:

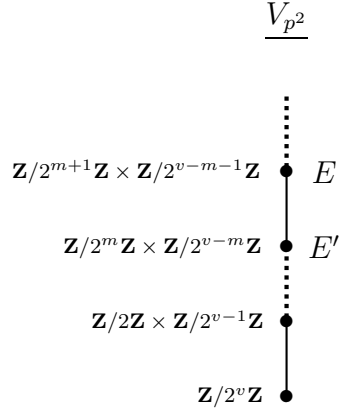
\begin{center}
\begin{tikzpicture}[scale=1.0,sizefont/.style={scale = 2}]
\draw[ultra thick] (1,5) node {};
\draw[ultra thick, dotted] (1,5) -- (1,4);
\draw[ultra thick] (1,4) node {${\bullet}$};
\draw[thick] (1,3) -- (1,4);
\draw[ultra thick] (1,3) node {${\bullet}$};
\draw[ultra thick] (1,2) node {${\bullet}$};
\draw[thick] (1,2) -- (1,1);
\draw[ultra thick, dotted] (1,1) node {${\bullet}$};
\draw[ultra thick, dotted] (1,3) -- (1,2);
\draw(-1,4) node {\tiny{$\Z/2^{m+1} \Z \times \Z/2^{v-m-1}\Z$}};
\draw(-0.65,3) node {\tiny{$\Z/2^m\Z \times \Z/2^{v-m}\Z$}};
\draw(-0.45,2) node {\tiny{$\Z/2\Z \times \Z/2^{v-1}\Z$}};
\draw(0.25,1) node {\tiny{$\Z/2^v\Z$}};
\draw(1.5,4) node {$E$};
\draw(1.5,3) node {$E'$};
\draw(1,6) node {\underline{$V_{p^2}$}};
\end{tikzpicture}
\end{center}

It follows that $E'(\F_{p^2})$ has full $2^m$-torsion, but not full $2^{m+1}$-torsion.  Therefore, $p$ is anomalous of defect $(m+1,m)$.
\end{proof}

In our final proof below we will employ this proposition several times by exhibiting a Frobenius matrix $F$ that sets the height of $V_p$ to be $m$ and situates $E$ on the crater of $V_p$.  It will then follow that $p$ is anomalous of defect $(m+1,m)$. We now proceed to the proof.

\subsection{The Proof of Theorem \ref{cm_thm}} \label{final_proof} 

\begin{proof}[Proof of Theorem \ref{cm_thm}]
We have reduced the proof of the theorem to the four cases outlined in Section \ref{1/12}.  In each case, switching the roles of $E$ and $E'$ if necessary, we choose $E$ and $E'$ so that the $\Q$-rational isogeny $E \to E'$ reduces to a descending isogeny modulo every good prime $p$.  Thus $c_{2^m}' = 0$ for all $m \geq 2$.  By Lemma \ref{CMkernel} we therefore have
\[
\PPP(E,E') = \frac{1}{8} \sum_{m =2}^\infty \frac{c_{2^m}}{4^{m-2}}.
\]
It remains to prove that $c_{2^m} = 1/2$ for all $m$.  

To do this, observe that $E$ and $E'$ satisfy the hypotheses of Proposition \ref{c2mprop}; in particular, both $G$ and $G'$ contain $-I$, for each group in Lemma \ref{CMkernel}.  Therefore, since $c_{2^m}'=0$, if we can show that $c_{2^m} \ne 0$, then it must be the case that $c_{2^m} = 1/2$.  
Proposition \ref{literal_zero} implies that 
it suffices to prove that there exist primes of defect $(m+1,m)$ for all $m\ge 2$. We  do this by exhibiting a single matrix in $G(2^m)$ representing $F$ at an anomalous prime. 
This will imply that 
\[
\PPP(E,E') = \frac{1}{16} \sum_{m \geq 2} \frac{1}{4^{m-2}} = \frac{1}{12},
\]
which will complete the proof.  We divide the remainder of the proof into cases according to the possible Galois images of Lemma \ref{CMkernel}.  However, in each case our strategy is identical.  Using the explicit matrix descriptions given in \cite{alvaro}, we show that, for every $m\geq 2$, each possible Galois image contains a Frobenius class satisfying the hypotheses of Proposition \ref{key_height_prop}.  Now we consider each scenario individually (grouping the two cases of $j$-invariant 0 into one case); for the remainder of the proof we fix $m \geq 2$.  We also import some notation from \cite{alvaro} and give brief descriptions when necessary.

\medskip

\noindent \fbox{\textbf{Case 1: $G = G' = N_{\delta,\phi}(2^\infty)$}} Recall from \cite[Thm.~1.1]{alvaro} that
\[
N_{\delta,\phi}(M) = \left\< \begin{pmatrix} -1 & 0 \\ \phi & 1\end{pmatrix}, C_{\delta,\phi}(M) \right\>,
\]
where 
\[
C_{\delta,\phi}(M) = \left\{ \begin{pmatrix} a + b\phi & b \\ \delta b & a \end{pmatrix} ~:~ a,b\in \Z/M\Z \text{ and } a^2 + ab \delta -\delta b^2 \in (\Z/M\Z)^\times \right\}.
\]
The quantities $\delta,\phi$ are determined by the elliptic curve $E$.  Suppose $E$ has CM by the order $\mathcal{O}_{K,f}$ of conductor $f$ in the imaginary quadratic field $K$; write $\Delta_Kf^2$ for the discriminant of $\mathcal{O}_{K,f}$.  Then the quantities $\delta,\phi, \Delta_Kf^2$, and $f$ are related by \cite[Thm.~1.1]{alvaro}:
\begin{enumerate}
\item \label{alv1} If $\Delta_Kf^2 \equiv 0 \pmod{4}$, or $M$ is odd, then $\delta = \Delta_Kf^2/4$ and $\phi = 0$.
\item \label{alv2} If $\Delta_Kf^2 \equiv 1 \pmod{4}$ and $M$ is even, then $\delta = (\Delta_K-1)f^2/4$ and $\phi = f$. 
\end{enumerate}
Finally, we have $N_{\delta,\phi}(2^\infty) = \varprojlim_{n} N_{\delta,\phi}(2^n)$. In $C_{\delta,\phi}(2^{m+1})$, consider the element 
\[
F = \begin{pmatrix} a + b\phi & b \\ \delta b & a \end{pmatrix} = \begin{pmatrix} -1 + 2^m\phi & 2^m \\ \delta 2^m & -1 \end{pmatrix} = -I + 2^m \begin{pmatrix} \phi & 1 \\ \delta  & 0 \end{pmatrix}.
\]
Equation \eqref{height_matrix} implies that the 
height of $V_p$ is $m + v_2(\phi^2 + 4\delta)$.   Since $F \equiv -I \pmod{2^m}$ then the level of $E$ is at least $m$.  Thus, if we can show that the height of $V_p$ is exactly $m$, then we may conclude that the level of $E$ is $m$.

Since $j(E) \in \Q$, we have that $\End(E)$ has class number 1 and, because $E$ lies in the $R_4(14)$ graph, $\OK \otimes \Q \simeq \Q(\sqrt{-7})$.  Thus $\Delta_K \equiv 1 \pmod{8}$.  Because $M$ is even, this puts us in case (\ref{alv2}) of the previous paragraph:
\[
\Delta_Kf^2 \equiv 1 \pmod{4}, \text{ and } \delta = (\Delta_K-1)f^2/4, \text{ and }  \phi = f.
\]
Thus, $\phi^2 + 4\delta = f^2 + (\Delta_K^2 - 1)$. But $\Delta_K^2-1 \equiv 0 \pmod{4}$.  Since $\Delta_K^2f^2 \equiv 1 \pmod{4}$, it must be the case that $\phi^2  = f^2 \equiv  1 \pmod{4}$.  Thus $v_2(\phi^2 + 4\delta) = 0$, so $h(V_p)=m$.  Therefore, the level of $E$ is exactly $m$ and Proposition \ref{key_height_prop} implies that $p$ has defect $(m+1,m)$.

\medskip

\noindent \textbf{\fbox{Case 2: $j(E) =0$}} Of the groups occurring in Lemma \ref{CMkernel}, two of them correspond to curves with $j$-invariant 0:  $G = \langle C_{-1,1}(2^\infty)^3,c_{-1}'\rangle$ and $G= \langle C_{-1,1}(2^\infty)^3,c_{1}'\rangle$.  The group $C_{-1,1}(2^\infty)^3$ contains all invertible matrices of the form
\[
\begin{pmatrix} a + b & b \\ -b & a \end{pmatrix}
\]
with $b$ even.  In particular, it contains the matrix
\[
-I + 2^m \begin{pmatrix} 1 & 1 \\ -1 & 0\end{pmatrix}.
\]
If such a matrix is in the conjugacy class of Frobenius for an elliptic curve $E$, then Equation \eqref{height_matrix} implies that the height of $V_p$ is
\[
m + v_2(\phi^2 + 4\delta) = m + v_2(1^2 -4) = m.
\]
Since $F\equiv -I\pmod{2^m}$, the level of $E$ is at least $m$, and we conclude that the level is exactly $m$.  Proposition \ref{key_height_prop} then implies that $p$ is anomalous of defect $(m+1,m)$. 

\medskip

\noindent \fbox{\textbf{Case 3:} $G_{2,a}$}  By the explicit description in \cite[Thm.~1.7]{alvaro}, the group $G_{2,a}$ contains the matrix 
\[
A = \begin{pmatrix} 1&2 \\ -2 & 1 \end{pmatrix}.
\]
It is also the case that elliptic curves with 2-adic image $G_{2,a}$ have $j$-invariant 1728.  

One computes that for $m \geq 1$, 
\[
A^{2^{m-1}}  = -I + 2^m \begin{pmatrix} v & u \\ -u & v \end{pmatrix}
\]
for odd integers $u,v$.  Because $j(E) = 1728$, we have that $\disc \mathcal{O}_0 = -4$.  Now suppose $p$ is a good prime with $F = A^{2^{m-1}}$.  Equation \eqref{height_matrix} implies that the height of $V_p$ is
\[
m + \frac{1}{2} v_2(((v-v)^2 - 4u^2)/(-4)) = m +  v_2(u^2)/2 = m.
\]
As in the previous cases, the level of $E$ must be $m$ and applying Proposition \ref{key_height_prop} completes the proof.
\end{proof}

\end{document}